\let\ORIlabel\label
\let\ORIrefstepcounter\refstepcounter
\AddToHook{package/hyperref/before}{%
  \let\label\ORIlabel
  \let\refstepcounter\ORIrefstepcounter}
\documentclass[onefignum,onetabnum]{siamart171218}

\usepackage{amssymb}
\usepackage{booktabs}
\usepackage{array}
\usepackage{placeins}

\newsiamthm{assumption}{Assumption}
\newsiamremark{example}{Example}
\newsiamremark{remark}{Remark}

\crefname{assumption}{Assumption}{Assumptions}
\Crefname{assumption}{Assumption}{Assumptions}

\newcommand{\R}{\mathbb{R}}
\newcommand{\E}{\mathbb{E}}
\newcommand{\PP}{\mathbb{P}}
\newcommand{\HH}{\mathcal{H}}
\newcommand{\F}{\mathcal{F}}
\newcommand{\inner}[2]{\langle #1,#2\rangle}
\newcommand{\norm}[1]{\lVert #1\rVert}
\DeclareMathOperator*{\argmin}{arg\,min}
\DeclareMathOperator{\Fix}{Fix}
\newcommand{\En}{\mathbb{E}_n}

\headers{Convergence of HalpernSGD}{Vittorio Colao and Katherine Rossella Foglia}
\title{On the Convergence of HalpernSGD%
\thanks{The authors contributed equally to this work.}}

\author{Vittorio Colao%
\thanks{Corresponding author. Department of Mathematics and Computer Science,
University of Calabria, Ponte P. Bucci, 30B, 87036 Arcavacata di Rende (CS),
Italy (\email{vittorio.colao@unical.it}).}
\and Katherine Rossella Foglia%
\thanks{Department of Mathematics and Computer Science, University of
Calabria, Ponte P. Bucci, 30B, 87036 Arcavacata di Rende (CS), Italy
(\email{katherine.foglia@unical.it}).}}

\begin{document}

\maketitle

\begin{abstract}
We study HalpernSGD for minimizing a convex Fr\'echet differentiable function
with Lipschitz gradient on a real Hilbert space. Its deterministic orbit
converges strongly to the anchor-selected minimizer \(P_S(u)\) under vanishing
steps that need not be square summable. For stochastic oracles, weighted
summability of the scaled martingale fluctuation and the \(\F_n\)-measurable bias yields
almost sure strong convergence, while expected last-iterate bounds require
suffix regularity of the error budget. In the unbiased case, for the critical
choice \(\varepsilon_n=n^{-1/2}\), \(\alpha_n=(\log n)/n\) eventually, a
uniform conditional moment bound of order \(p>4\) gives both conclusions. We
also discuss finite-difference and zeroth-order applications for which scaled,
rather than raw, noise is the natural quantity.
\end{abstract}




\section{Introduction}\label{sec:intro}

Let $\HH$ be a real Hilbert space and let $f\colon\HH\to\R$ be convex and
Fr\'echet differentiable, with $L$-Lipschitz gradient and nonempty solution set
\(S:=\argmin f\). Given an anchor $u\in\HH$, Halpern schemes select the
distinguished minimizer \(x^*:=P_S(u)\) rather than an arbitrary point of
$S$~\cite{halpern1967,wittmann1992,xu2002another}. We consider
\begin{equation}\label{eq:intro-det}
  x_{n+1}
  =
  \alpha_n u + (1-\alpha_n)\bigl(x_n-\varepsilon_n\nabla f(x_n)\bigr),
\end{equation}
where $(\alpha_n)\subset(0,1)$ and $(\varepsilon_n)\subset(0,2/L]$. Since
$I-\varepsilon\nabla f$ is nonexpansive and has fixed-point set $S$ for
$0<\varepsilon\le2/L$~\cite{baillon1977,bauschke2010baillon},
\eqref{eq:intro-det} is a Halpern iteration with varying maps. Its stochastic
counterpart is
\begin{equation}\label{eq:intro-stoch}
  X_{n+1}
  =
  \alpha_n u + (1-\alpha_n)\bigl(X_n-\varepsilon_n G_n\bigr),
  \qquad
  G_n=\nabla f(X_n)+U_n,
\end{equation}
where $U_n$ is the oracle error. We call \eqref{eq:intro-stoch} HalpernSGD; its
name and empirical neural-network variants originated
in~\cite{foglia2024halpernsgd,colao2025optimizer}.

Under \eqref{eq:alpha-base}--\eqref{eq:coupling}, the deterministic orbit
converges strongly to $P_S(u)$ without assuming
$\sum_n\varepsilon_n^2<\infty$; the proof combines an anchor-sensitive block
estimate with a Maing\'e index argument. A pathwise comparison then gives almost sure
convergence of HalpernSGD under weighted summability of the scaled martingale
fluctuation and the \(\F_n\)-measurable bias. Expected last-iterate estimates instead
require suffix regularity, and the two conclusions hold simultaneously in the
unbiased critical regime \(\varepsilon_n=n^{-1/2}\),
\(\alpha_n=(\log n)/n\) eventually, under a uniform conditional $p$-moment
bound with $p>4$. Finite-difference and zeroth-order examples, including
Kiefer--Wolfowitz and SPSA  
under the additional conditions stated in \cref{subsec:random-direction-application}, illustrate why the scaled noise
can be summable even when the raw variance is
not~\cite{kiefer1952stochastic,spall1992multivariate}.

These results complement recent work on deterministic and stochastic Halpern
schemes~\cite{diakonikolas2020halpern,park2022exact,contreras2023optimal,
bravo2026minimax,alacaoglu2025weakervar,bravo2024stochastic,
cai2022stochastic,bravocominetti2024,pischke2024asymptotic}; see
\cref{sec:comparison} for a detailed comparison.

\section{Preliminaries}\label{sec:prelim}

Throughout, $\HH$ is a real Hilbert space with inner product
$\inner{\cdot}{\cdot}$ and norm $\norm{\cdot}$, and $I$ is the identity. We use
standard asymptotic notation and the convex-analytic conventions of
\cite{bauschke2011convex}. For stochastic results, $\HH$ is additionally
separable and \((\Omega,\F,(\F_n)_{n\ge0},\PP)\) is a complete filtered
probability space. Hilbert-valued random variables are strongly measurable,
and expectations and conditional expectations are Bochner ones; we write
\(\E_n[\cdot]:=\E[\cdot\mid\F_n]\).

We consider the convex minimization problem
\begin{equation}\label{eq:prob}
  \min_{x\in\HH} f(x),
\end{equation}
where $f\colon\HH\to\R$ is convex and Fr\'echet differentiable with
$L$-Lipschitz gradient, $L>0$,
\begin{equation}\label{eq:Lsmooth}
  \norm{\nabla f(x)-\nabla f(y)}\le L\norm{x-y},
  \qquad x,y\in\HH,
\end{equation}
and the solution set
\(S:=\argmin_{x\in\HH} f(x)\)
is nonempty. Since $S$ is closed and convex, for every $u\in\HH$ the projection
\begin{equation}\label{eq:xstar}
  x^*:=P_S(u)
\end{equation}
is well defined and satisfies~\cite{bauschke2011convex}
\[
  \inner{u-x^*}{x-x^*}\le 0,
  \qquad x\in S.
\]

By the Baillon--Haddad theorem~\cite{baillon1977,bauschke2010baillon}, the
gradient is $1/L$-cocoercive:
\begin{equation}\label{eq:cocoercive}
  \inner{\nabla f(x)-\nabla f(y)}{x-y}
  \ge
  \tfrac1L\norm{\nabla f(x)-\nabla f(y)}^2,
  \qquad x,y\in\HH.
\end{equation}
Consequently, $T_\varepsilon:=I-\varepsilon\nabla f$ is nonexpansive for
$0<\varepsilon\le2/L$, and
\begin{equation}\label{eq:Fix}
  \Fix(T_\varepsilon)=S,
  \qquad 0<\varepsilon\le 2/L.
\end{equation}
We also use the standard descent and convexity inequalities
\begin{align*}
  f(y)&\le f(x)+\inner{\nabla f(x)}{y-x}+\tfrac L2\norm{y-x}^2,
        &&x,y\in\HH,\\
  f(x)-f(z)&\le\inner{\nabla f(x)}{x-z},
        &&x\in\HH,\ z\in S.
\end{align*}

Let $(\alpha_n)\subset(0,1)$ and $(\varepsilon_n)\subset(0,2/L]$. Our baseline
conditions are
\begin{equation}\label{eq:alpha-base}
  \alpha_n\to0,
  \qquad
  \sum_{n=0}^{\infty}\alpha_n=+\infty,
\end{equation}
the vanishing of the gradient step,
\begin{equation}\label{eq:eps-base}
  \varepsilon_n\to0,
\end{equation}
and the scale separation
\begin{equation}\label{eq:coupling}
  \frac{\varepsilon_n}{\alpha_n}\to+\infty.
\end{equation}
Together with \(\sum_n\alpha_n=+\infty\), \eqref{eq:coupling} implies
\(\sum_n\varepsilon_n=+\infty\). Further square-summability or suffix
conditions are imposed only where needed.

For \eqref{eq:intro-stoch}, assume
\(X_0\in L^2(\Omega,\F_0;\HH)\) and
\(U_n\in L^2(\Omega,\F_{n+1};\HH)\); then $(X_n)$ is adapted. We use the
decomposition
\begin{equation}
\label{eq:noise-decomposition}
  U_n=b_n+\xi_n,
  \qquad
  b_n:=\E[U_n\mid\F_n],
  \qquad
  \E[\xi_n\mid\F_n]=0,
\end{equation}
where $b_n$ is $\F_n$-measurable but generally random (and may depend on $X_n$),
while $\xi_n$ is conditionally centered. Define
\begin{equation}
\label{eq:def-vn}
  v_n:=\varepsilon_n^2\,\E\bigl[\norm{\xi_n}^2\mid\F_n\bigr],
\end{equation}
and
\begin{equation}
\label{eq:def-bias-remainder}
  \beta_n
  :=
  \varepsilon_n
  \bigl(\E\norm{b_n}^2\bigr)^{1/2}
  +
  \varepsilon_n^2
  \E\norm{b_n}^2.
\end{equation}
In the unbiased case $b_n=0$, $\xi_n=U_n$, $\beta_n=0$, and
$v_n=\varepsilon_n^2\E_n\norm{U_n}^2$.

\section{Deterministic Halpern--gradient dynamics}
\label{sec:det}

For \eqref{eq:intro-det}, set $x^*=P_S(u)$ and
\(A_n:=\norm{x_n-x^*}^2\),
\(Q_n:=\inner{\nabla f(x_n)}{x_n-x^*}\), and
\(P_n:=\inner{u-x^*}{x_n-x^*}\). Convexity and cocoercivity give
\begin{equation}
\label{eq:gap-by-Qn-det}
        0\le f(x_n)-f(x^*)\le Q_n.
\end{equation}
\begin{equation}
\label{eq:cocoercivity-det}
        \norm{\nabla f(x_n)}^2\le LQ_n.
\end{equation}

\begin{lemma}
\label[lemma]{lem:basic-dissipative}
Assume that
\(0<\alpha_n<1\) and \(0<\varepsilon_n\le\frac2L.\)
Then the sequence generated by~\eqref{eq:intro-det} satisfies
\begin{equation}
\label{eq:basic-dissipative}
        A_{n+1}
        +(1-\alpha_n)\varepsilon_n(2-L\varepsilon_n)Q_n
        \le
        (1-\alpha_n)A_n+\alpha_n\norm{u-x^*}^2.
\end{equation}
Consequently,
\(A_n\le \max\{A_0,\norm{u-x^*}^2\}\) for $n\ge0$.
In particular, \((x_n)\) and \((\nabla f(x_n))\) are bounded.
\end{lemma}
\begin{proof}
Set \(y_n:=x_n-\varepsilon_n\nabla f(x_n)\). Since
\(x_{n+1}=\alpha_nu+(1-\alpha_n)y_n\), convexity of $\norm{\cdot}^2$ and
\eqref{eq:cocoercivity-det} yield
\[
 A_{n+1}\le\alpha_n\norm{u-x^*}^2+(1-\alpha_n)\norm{y_n-x^*}^2,
 \quad
 \norm{y_n-x^*}^2\le A_n-\varepsilon_n(2-L\varepsilon_n)Q_n.
\]
This proves \eqref{eq:basic-dissipative}.
After dropping its nonnegative
term, an induction yields the stated bound on \(A_n\);
Lipschitz continuity then bounds $\nabla f(x_n)$.
\end{proof}

Under \eqref{eq:alpha-base}--\eqref{eq:coupling}, set
\(h_n := (1-\alpha_n)^2\varepsilon_n(2-L\varepsilon_n),\)
and
\begin{equation}
\label{eq:def-anchor-block-det}
        \mathcal B_n
        :=
        A_n
        -
        2(1-\alpha_n)P_n
        +
        \frac{h_n}{\alpha_n}Q_n.
\end{equation}

\begin{lemma}
\label[lemma]{lem:anchor-block-det}
Assume~\eqref{eq:alpha-base}--\eqref{eq:coupling}. Then there exists a
nonnegative sequence \((c_n)\) such that
\(\frac{c_n}{\alpha_n}\to0\)
and
\begin{equation}
\label{eq:anchor-block-det}
        A_{n+1}+\alpha_n\mathcal B_n
        \le
        A_n+c_n.
\end{equation}
Moreover, let \(M:\mathbb N\to\mathbb N\) be any map such that
\(M(N)\ge N\), and define
\(\Gamma_N:=\sum_{n=N}^{M(N)}\alpha_n.\)
If
\(\Gamma_N\longrightarrow+\infty\), as $N\to\infty$,
then
\begin{equation}
\label{eq:min-anchor-block-det}
        \min_{N\le n\le M(N)}
        \mathcal B_n
        \le
        \frac{
            A_N+\sum_{n=N}^{M(N)}c_n
        }{
            \sum_{n=N}^{M(N)}\alpha_n
        }.
\end{equation}
The right-hand side converges to \(0\). Hence
\[
        \limsup_{N\to\infty}
        \min_{N\le n\le M(N)}
        \mathcal B_n
        \le0 .
\]
In particular, the assumption \(\sum_n\alpha_n=+\infty\) guarantees the
existence of admissible blocks: for each \(N\), choose \(M(N)\ge N\) so that
\(\sum_{n=N}^{M(N)}\alpha_n\ge N\).
Then \(\Gamma_N\ge N\to+\infty\), and hence
\begin{equation}
\label{eq:liminf-anchor-block-det}
        \liminf_{n\to\infty}\mathcal B_n\le0.
\end{equation}
\end{lemma}

\begin{proof}
From~\eqref{eq:intro-det},
\[
        x_{n+1}-x^*
        =
        \alpha_n(u-x^*)
        +
        (1-\alpha_n)
        \bigl(x_n-x^*-\varepsilon_n\nabla f(x_n)\bigr).
\]
Expanding the square,
\[
\begin{aligned}
        A_{n+1}
        &=
        \alpha_n^2\norm{u-x^*}^2
        +
        (1-\alpha_n)^2
        \norm{x_n-x^*-\varepsilon_n\nabla f(x_n)}^2              \\
        &\quad
        +
        2\alpha_n(1-\alpha_n)
        \inner{u-x^*}{x_n-x^*-\varepsilon_n\nabla f(x_n)} .
\end{aligned}
\]
The middle term satisfies, by~\eqref{eq:cocoercivity-det},
\[
\begin{aligned}
        \norm{x_n-x^*-\varepsilon_n\nabla f(x_n)}^2
        &=
        A_n
        -
        2\varepsilon_nQ_n
        +
        \varepsilon_n^2\norm{\nabla f(x_n)}^2                    \\
        &\le
        A_n-\varepsilon_n(2-L\varepsilon_n)Q_n.
\end{aligned}
\]
Therefore
\[
\begin{aligned}
        A_{n+1}
        &\le
        \alpha_n^2\norm{u-x^*}^2
        +
        (1-\alpha_n)^2A_n
        -
        (1-\alpha_n)^2\varepsilon_n(2-L\varepsilon_n)Q_n         \\
        &\quad
        +
        2\alpha_n(1-\alpha_n)P_n
        -
        2\alpha_n(1-\alpha_n)\varepsilon_n
        \inner{u-x^*}{\nabla f(x_n)} .
\end{aligned}
\]
By \cref{lem:basic-dissipative}, \((\nabla f(x_n))\) is bounded. Hence there
exists \(R_*>0\) such that, for $n\ge0$,
\(\bigl| \inner{u-x^*}{\nabla f(x_n)} \bigr| \le R_*\). Thus
\begin{equation}
\label{eq:diseqAm}
        A_{n+1}
        \le
        (1-\alpha_n)^2A_n
        +
        2\alpha_n(1-\alpha_n)P_n
        -
        h_nQ_n
        +
        c_n,
\end{equation}
where
\(c_n := \alpha_n^2\norm{u-x^*}^2 + 2\alpha_n(1-\alpha_n)\varepsilon_n R_*\).
In particular,
\[
        c_n=O(\alpha_n^2+\alpha_n\varepsilon_n),
        \qquad
        \frac{c_n}{\alpha_n}
        =
        O(\alpha_n+\varepsilon_n)\to0.
\]
Since
\((1-\alpha_n)^2 = 1-2\alpha_n+\alpha_n^2 \le 1-\alpha_n,\)
 we obtain
\[
\begin{aligned}
        A_{n+1}
        &\le
        (1-\alpha_n)A_n
        +
        2\alpha_n(1-\alpha_n)P_n
        -
        h_nQ_n
        +
        c_n                                                     \\
        &=
        A_n
        -
        \alpha_n
        \left[
            A_n
            -
            2(1-\alpha_n)P_n
            +
            \frac{h_n}{\alpha_n}Q_n
        \right]
        +
        c_n                                                    =
        A_n-\alpha_n\mathcal B_n+c_n.
\end{aligned}
\]
This proves~\eqref{eq:anchor-block-det}.
We now apply \cref{lem:tail-weighted-telescoping} with
$\varrho_n=1$, $w_n=\alpha_n$, $b_n=\mathcal B_n$, and $r_n=c_n$.
For any \(M(N)\ge N\), it gives \eqref{eq:min-anchor-block-det}.
If \(\Gamma_N:=\sum_{n=N}^{M(N)}\alpha_n\to+\infty\), then
\((A_N)\) is bounded by \cref{lem:basic-dissipative}, so
\(A_N/\Gamma_N\to0\), and
\[
        \frac{
            \sum_{n=N}^{M(N)}c_n
        }{
            \sum_{n=N}^{M(N)}\alpha_n
        }
        \le
        \sup_{n\ge N}\frac{c_n}{\alpha_n}
        \to0.
\]
Thus the right-hand side of \eqref{eq:min-anchor-block-det} tends to zero,
and the asymptotic assertion of \cref{lem:tail-weighted-telescoping} gives the
claimed limsup estimate.
Moreover, \eqref{eq:liminf-anchor-block-det} follows by choosing, for each \(N\), an index \(M(N)\ge N\) such that
\(\sum_{n=N}^{M(N)}\alpha_n\ge N\).
Such a choice is possible because \(\sum_n\alpha_n=+\infty\), and the selected
indices belong to intervals whose left endpoint tends to infinity.
\end{proof}

The convergence proof uses the following lemma of Maingé~\cite{mainge2008}.

\begin{lemma}
\label[lemma]{lem:mainge}
Let \((s_n)\) be a real sequence which is not eventually nonincreasing, and
choose an index \(n_0\) such that \(s_{n_0}<s_{n_0+1}\). For \(n\ge n_0\),
define
\[
  \tau(n):=\max\{k\in\{n_0,\ldots,n\}:\ s_k<s_{k+1}\}.
\]
Then
\(\tau(n)\to+\infty\), \(s_{\tau(n)}<s_{\tau(n)+1},\)
and
\(s_n\le s_{\tau(n)+1}\)
for all sufficiently large \(n\).
\end{lemma}

\begin{theorem}
\label{thm:det-strong-convergence}
Let \((x_n)\) be generated by~\eqref{eq:intro-det}. If
\eqref{eq:alpha-base}--\eqref{eq:coupling} hold, then
\(x_n\to x^*=P_S(u)\) strongly in \(\HH\).
\end{theorem}

\begin{proof}
Recall that
\(A_n=\norm{x_n-x^*}^2.\)
We distinguish two cases.

\smallskip
\noindent
Case 1: \((A_n)\) is eventually nonincreasing.
Then
\(A_n\to\ell \) for some \(\ell\ge0.\)
By \cref{lem:anchor-block-det},
\(\liminf_{n\to\infty}\mathcal B_n\le0.\)
Since \(A_n\ge0\), \(Q_n\ge0\), and \((P_n)\) is bounded, the sequence
\((\mathcal B_n)\) is bounded from below. Hence
\(\eta:=\liminf_{n\to\infty}\mathcal B_n\)
is finite and satisfies \(\eta\le0\). We may therefore choose a subsequence
\((n_k)\) such that
\(\mathcal B_{n_k}\to\eta\).
In the decomposition
\(\mathcal B_n = A_n - 2(1-\alpha_n)P_n + \frac{h_n}{\alpha_n}Q_n,\)
the sequences \((A_n)\) and \((P_n)\) are bounded because \((x_n)\) is bounded.
Moreover,
\[
        \frac{h_n}{\alpha_n}
        =
        \frac{(1-\alpha_n)^2\varepsilon_n(2-L\varepsilon_n)}
             {\alpha_n}
        \sim
        2\frac{\varepsilon_n}{\alpha_n}
        \to+\infty .
\]
We claim that \(Q_{n_k}\to0\). Indeed, if this were false, then, along a
subsequence, \(Q_{n_k}\) would be bounded away from zero. Since
\(Q_n\ge0\) and \(h_n/\alpha_n\to+\infty\), the term
\(\frac{h_{n_k}}{\alpha_{n_k}}Q_{n_k}\)
would then diverge to \(+\infty\). The remaining terms in
\(\mathcal B_{n_k} = A_{n_k}-2(1-\alpha_{n_k})P_{n_k} + \frac{h_{n_k}}{\alpha_{n_k}}Q_{n_k}\)
are bounded, because \((A_n)\) and \((P_n)\) are bounded. Hence
\(\mathcal B_{n_k}\to+\infty\), contradicting the boundedness from above of
\((\mathcal B_{n_k})\).

By reflexivity of \(\HH\) and weak sequential compactness of bounded sets in
Hilbert spaces, we may pass to a weakly convergent subsequence, still denoted
by \((x_{n_k})\), such that
\[
        x_{n_k}\rightharpoonup p.
\]
By~\eqref{eq:gap-by-Qn-det},
\(0\le f(x_{n_k})-f(x^*)\le Q_{n_k}\to0.\)
Since \(f\) is convex and continuous, it is weakly lower semicontinuous.
Therefore
\(f(p)\le \liminf_{k\to\infty}f(x_{n_k})=f(x^*).\)
Hence \(p\in S\). By the projection characterization of \(x^*=P_S(u)\),
\(\inner{u-x^*}{p-x^*}\le0.\)
Moreover, since \(x_{n_k}\rightharpoonup p\), the linear functional
\(x\mapsto\inner{u-x^*}{x-x^*}\) gives
\(P_{n_k} = \inner{u-x^*}{x_{n_k}-x^*} \to \inner{u-x^*}{p-x^*} \le0.\)
Using the definition of \(\mathcal B_n\) and \(Q_{n_k}\ge0\),
\[
\begin{aligned}
        A_{n_k}
        &=
        \mathcal B_{n_k}
        +
        2(1-\alpha_{n_k})P_{n_k}
        -
        \frac{h_{n_k}}{\alpha_{n_k}}Q_{n_k}       \le
        \mathcal B_{n_k}
        +
        2(1-\alpha_{n_k})P_{n_k}.
\end{aligned}
\]
Taking the limit and using \(A_{n_k}\to\ell\), we obtain
\(\ell \le \eta+2\inner{u-x^*}{p-x^*} \le0.\)
Since \(\ell\ge0\), it follows that \(\ell=0\). Therefore \(A_n\to0\).

\smallskip
\noindent
Case 2: \((A_n)\) is not eventually nonincreasing.
Choose \(n_0\) such that \(A_{n_0}<A_{n_0+1}\), and let
\(\tau(n)\), \(n\ge n_0\), be the Maingé index of \((A_n)\) defined in
\cref{lem:mainge}. Then
\(\tau(n)\to+\infty, \qquad A_{\tau(n)}<A_{\tau(n)+1},\)
and
\(A_n\le A_{\tau(n)+1}\)
for all large \(n\).

Set \(m:=\tau(n)\). From~\eqref{eq:basic-dissipative} and
\(A_m<A_{m+1}\), as $h_m\leq (1-\alpha_m)\varepsilon_m(2-L\varepsilon_m)$, we obtain
\[
\begin{aligned}
          A_m+h_mQ_m
        & \le
        A_{m+1}+h_mQ_m
        \le
        A_{m+1}+(1-\alpha_m)\varepsilon_m(2-L\varepsilon_m)Q_m
        \\
        &\le
        (1-\alpha_m)A_m+\alpha_m\norm{u-x^*}^2.
\end{aligned}
\]
Hence
\(\alpha_mA_m+h_mQ_m \le \alpha_m\norm{u-x^*}^2.\)
In particular,
\(Q_m \le \frac{\alpha_m}{h_m}\norm{u-x^*}^2.\)
Since
\(h_m = (1-\alpha_m)^2\varepsilon_m(2-L\varepsilon_m) \sim 2\varepsilon_m\)
and \(\alpha_m/\varepsilon_m\to0\), we get
\(Q_{\tau(n)}\to0.\)

Set \(r:=\limsup_{n\to\infty}P_{\tau(n)}\), which is finite because
\((x_n)\) is bounded. Choose a subsequence \((n_j)\) such that
\(P_{\tau(n_j)}\to r\). The sequence \((x_{\tau(n_j)})\) is bounded, so, after
passing to a further subsequence, there is \(p\in\HH\) such that
\[
        x_{\tau(n_j)}\rightharpoonup p.
\]
By~\eqref{eq:gap-by-Qn-det},
\(0\le f(x_{\tau(n_j)})-f(x^*)\le Q_{\tau(n_j)}\to0.\)
Weak lower semicontinuity gives \(f(p)\le f(x^*)\), hence \(p\in S\). Thus
\(\inner{u-x^*}{p-x^*}\le0.\)
By weak convergence, \(r\le0\), and hence
\begin{equation}
\label{eq:limsup-anchor-mainge}
        \limsup_{n\to\infty}
        \inner{u-x^*}{x_{\tau(n)}-x^*}
        \le0.
\end{equation}

Using the direct expansion from the proof of
\cref{lem:anchor-block-det} (see \eqref{eq:diseqAm}), for \(m=\tau(n)\),
\(A_{m+1} \le (1-\alpha_m)^2A_m + 2\alpha_m(1-\alpha_m)P_m - h_mQ_m + c_m.\)
Since \(A_m<A_{m+1}\) and \(-h_mQ_m\le0\),
\(\bigl[1-(1-\alpha_m)^2\bigr]A_m \le 2\alpha_m(1-\alpha_m)P_m+c_m.\)
As
\(1-(1-\alpha_m)^2 = \alpha_m(2-\alpha_m),\)
we have
\(A_m \le \frac{2(1-\alpha_m)}{2-\alpha_m}P_m + \frac{c_m}{\alpha_m(2-\alpha_m)}.\)
Since \(c_m/\alpha_m\to0\),
\(A_m \le \frac{2(1-\alpha_m)}{2-\alpha_m}P_m+o(1).\)
Taking the limsup along \(m=\tau(n)\) and using
\eqref{eq:limsup-anchor-mainge}, we get
\(\limsup_{n\to\infty}A_{\tau(n)}\le0.\)
Since \(A_{\tau(n)}\ge0\), this implies
\(A_{\tau(n)}\to0.\)
By~\eqref{eq:basic-dissipative},
\(A_{\tau(n)+1} \le (1-\alpha_{\tau(n)})A_{\tau(n)} + \alpha_{\tau(n)}\norm{u-x^*}^2 \to0.\)
Finally, by \cref{lem:mainge},
\(A_n\le A_{\tau(n)+1}\)
for all sufficiently large \(n\). Hence \(A_n\to0\).

In both cases,
\(A_n=\norm{x_n-x^*}^2\to0.\)
Therefore
\(x_n\to x^* \) strongly in $\HH.$
\end{proof}

We now derive a deterministic last-iterate estimate for the functional gap.
The statement is given in a suffix-regular form because the stochastic
last-iterate estimate will have the same structure, with the additional
variance contribution.

We follow the standard convention and write
\(A_K=\widetilde O(B_K)\)
to mean that, for some \(r\ge0\),
\(A_K=O\bigl(B_K(\log K)^r\bigr)\).

\begin{assumption}
\label[assumption]{ass:det-suffix-rate-regime}
The sequence \((\varepsilon_n)\) is eventually nonincreasing, and
\begin{align}
        \frac{
            1+\sum_{k=1}^{K}(\alpha_k+\varepsilon_k^2)
        }{K}
        &=
        \widetilde O(\varepsilon_K^2),
        \label{eq:rate-global-det}
        \\
        \sum_{\ell=1}^{K-1}
        \frac1{\ell(\ell+1)}
        \sum_{k=K-\ell}^{K}
        (\alpha_k+\varepsilon_k^2)
        &=
        \widetilde O(\varepsilon_K^2).
        \label{eq:rate-suffix-det}
\end{align}
\end{assumption}

\begin{theorem}
\label{thm:det-last-iterate-general}
Suppose that the hypotheses of \cref{thm:det-strong-convergence} hold and that
\cref{ass:det-suffix-rate-regime} is satisfied. Then
\(f(x_K)-f(x^*)=\widetilde O(\varepsilon_K)\).
\end{theorem}

\begin{proof}
Set
\(\Delta_k:=f(x_k)-f(x^*).\)
We first establish a one-step inequality with an arbitrary comparator
\(z\in\HH\). Let
\(y_k:=x_k-\varepsilon_k\nabla f(x_k).\)
By convexity of the squared norm,
\(\norm{x_{k+1}-z}^2 \le \alpha_k\norm{u-z}^2 + (1-\alpha_k)\norm{y_k-z}^2.\)
Moreover,
\[
        \norm{y_k-z}^2
        =
        \norm{x_k-z}^2
        -
        2\varepsilon_k\inner{\nabla f(x_k)}{x_k-z}
        +
        \varepsilon_k^2\norm{\nabla f(x_k)}^2.
\]
By convexity,
\(f(x_k)-f(z) \le \inner{\nabla f(x_k)}{x_k-z}.\)
Therefore
\begin{equation}
\label{eq:one-step-det-comparator}
\begin{aligned}
        2(1-\alpha_k)\varepsilon_k\bigl(f(x_k)-f(z)\bigr)
        &\le
        (1-\alpha_k)\norm{x_k-z}^2
        -
        \norm{x_{k+1}-z}^2
        +
        \alpha_k\norm{u-z}^2  \\
        &\quad
        +
        (1-\alpha_k)\varepsilon_k^2
        \norm{\nabla f(x_k)}^2
        \le
        \big( \norm{x_k-z}^2
        -
        \norm{x_{k+1}-z}^2 \big)    \\
        & \quad
        + \alpha_k\norm{u-z}^2
        + \varepsilon_k^2
        \norm{\nabla f(x_k)}^2.
\end{aligned}
\end{equation}
By \cref{lem:basic-dissipative}, \((x_k)\) and \((\nabla f(x_k))\) are
bounded. Let
\(G:=\sup_k\norm{\nabla f(x_k)}<+\infty.\)
Since \(\alpha_k\to0\), choose \(k_0\) such that
\(2(1-\alpha_k)\ge1\) for every \(k\ge k_0\).

Choose \(z=x^*\) and note that
\(\Delta_k=f(x_k)-f(x^*)\ge0\). Summing
\eqref{eq:one-step-det-comparator} from \(k=k_0\) to \(K\) and dropping the
final negative distance term gives
\[
        \sum_{k=k_0}^{K}2(1-\alpha_k)\varepsilon_k\Delta_k
        \le
        \norm{x_{k_0}-x^*}^2
        +
        \norm{u-x^*}^2\sum_{k=k_0}^{K}\alpha_k
        +
        G^2\sum_{k=k_0}^{K}\varepsilon_k^2.
\]
Adding the finitely many nonnegative terms with \(k<k_0\) and enlarging a
constant \(C>0\), we obtain
\[
        \sum_{k=1}^{K}\varepsilon_k\Delta_k
        \le
        C\left(
            1+\sum_{k=1}^{K}(\alpha_k+\varepsilon_k^2)
        \right).
\]
Dividing by \(K\),
\begin{equation}
\label{eq:average-det-rate}
        \frac1K
        \sum_{k=1}^{K}\varepsilon_k\Delta_k
        \le
        C
        \frac{
            1+\sum_{k=1}^{K}(\alpha_k+\varepsilon_k^2)
        }{K}
        =
        \widetilde O(\varepsilon_K^2),
\end{equation}
by~\eqref{eq:rate-global-det}.

For the suffix estimate, fix \(K\) and \(1\le\ell\le K-1\), and set
\(m=K-\ell\). Apply \eqref{eq:one-step-det-comparator} with \(z=x_m\) and sum
from \(k=m\) to \(K\). The distance terms give

\[
\sum_{k=m}^{K}
\bigl((1-\alpha_k)\|x_k-x_m\|^2-\|x_{k+1}-x_m\|^2\bigr)
\le
-\|x_{K+1}-x_m\|^2
\le0,
\]
where the term \(k=m\) vanishes because \(x_m-x_m=0\). Since \((x_n)\) is
bounded, the remaining anchor term satisfies \(\|u-x_m\|^2\le C\). Hence
\begin{equation}
\label{eq:quasi-suffix-det-estimate}
        \sum_{k=K-\ell}^{K}
        2(1-\alpha_k)\varepsilon_k
        \bigl(f(x_k)-f(x_{K-\ell})\bigr)
        \le
        C
        \sum_{k=K-\ell}^{K}
        (\alpha_k+\varepsilon_k^2).
\end{equation}

Since
\((x_n)\) is bounded and \(f\) is \(L\)-smooth, the values \(f(x_n)\) are
bounded; hence
\(\bigl|f(x_k)-f(x_m)\bigr|\le C\).
Dividing \eqref{eq:quasi-suffix-det-estimate} by two gives an upper bound for
the sum with coefficient \((1-\alpha_k)\varepsilon_k\). Moreover,
\[
\begin{aligned}
  \sum_{k=m}^{K}\varepsilon_k\bigl(f(x_k)-f(x_m)\bigr)
  &=
  \sum_{k=m}^{K}(1-\alpha_k)\varepsilon_k
       \bigl(f(x_k)-f(x_m)\bigr) \\
  &\quad+
  \sum_{k=m}^{K}\alpha_k\varepsilon_k
       \bigl(f(x_k)-f(x_m)\bigr).
\end{aligned}
\]
The second sum is bounded above in absolute value by
\(C\sum_{k=m}^{K}\alpha_k\varepsilon_k\), which is at most a constant
multiple of \(\sum_{k=m}^{K}\alpha_k\) because
\(\varepsilon_k\le2/L\). Hence
\[
        \sum_{k=m}^{K}
        \varepsilon_k\bigl(f(x_k)-f(x_m)\bigr)
        \le
        C\sum_{k=m}^{K}(\alpha_k+\varepsilon_k^2).
\]
Here, as throughout the proof, \(C>0\) denotes a generic constant, possibly enlarged from line to line, and independent of \(m\), \(K\), and \(\ell\). Thus, with \(m=K-\ell\),
\begin{equation}
\label{eq:suffix-det-estimate}
        \sum_{k=K-\ell}^{K}
        \varepsilon_k
        \bigl(f(x_k)-f(x_{K-\ell})\bigr)
        \le
        C
        \sum_{k=K-\ell}^{K}
        (\alpha_k+\varepsilon_k^2).
\end{equation}

Since
\(\Delta_k-\Delta_{K-\ell} = f(x_k)-f(x_{K-\ell}),\)
\eqref{eq:suffix-det-estimate} gives
\begin{equation} \label{eq:suffix-weighted-gap-det}
 \sum_{k=m}^{K}
        \varepsilon_k(\Delta_k-\Delta_m)
        \le
        C
        \sum_{k=m}^{K}
        (\alpha_k+\varepsilon_k^2).
\end{equation}

Choose \(N_0\) such that \((\varepsilon_k)_{k\ge N_0}\) is nonincreasing.
Since  \(\Delta_m\ge0\),  if \(m\ge N_0\) then
\[
\begin{aligned}
        \sum_{k=m}^{K}
        \bigl(\varepsilon_k\Delta_k-\varepsilon_m\Delta_m\bigr)
        &\le
        \sum_{k=m}^{K}
        \varepsilon_k(\Delta_k-\Delta_m).
\end{aligned}
\]
Therefore, if \(m=K-\ell\ge N_0\), applying
\eqref{eq:suffix-weighted-gap-det} gives
\[
        \sum_{k=K-\ell}^{K}
        \bigl(
        \varepsilon_k\Delta_k
        -
        \varepsilon_{K-\ell}\Delta_{K-\ell}
        \bigr)
        \le
        C
        \sum_{k=K-\ell}^{K}
        (\alpha_k+\varepsilon_k^2).
\]
Apply \cref{lem:suffix-averaging}, in the form of the identity
\eqref{eq:suffix-averaging-identity}, to the bounded sequence
\((\varepsilon_k\Delta_k)\). For the terms with
\(K-\ell\ge N_0\), use the preceding estimate. There are only finitely many
remaining left endpoints \(m=K-\ell<N_0\); because the sequence is bounded,
their total contribution to the suffix-averaging identity is \(O(1/K)\).
Consequently,
\[
\begin{aligned}
        \varepsilon_K\Delta_K
        &=
        \frac1K\sum_{k=1}^{K}\varepsilon_k\Delta_k
        +
        \sum_{\ell=1}^{K-1}
        \frac1{\ell(\ell+1)}
        \sum_{k=K-\ell}^{K}
        \bigl(
        \varepsilon_k\Delta_k
        -
        \varepsilon_{K-\ell}\Delta_{K-\ell}
        \bigr)                                                \\
        &\le
        \frac1K\sum_{k=1}^{K}\varepsilon_k\Delta_k
        +
        C
        \sum_{\ell=1}^{K-1}
        \frac1{\ell(\ell+1)}
        \sum_{k=K-\ell}^{K}
        (\alpha_k+\varepsilon_k^2)
        +O(1/K).
\end{aligned}
\]
The first term is \(\widetilde O(\varepsilon_K^2)\) by
\eqref{eq:average-det-rate}, and the second is
\(\widetilde O(\varepsilon_K^2)\) by
\eqref{eq:rate-suffix-det}. Finally, \eqref{eq:rate-global-det} implies
\(1/K=\widetilde O(\varepsilon_K^2)\), since the numerator on its left-hand
side is at least one. Hence
\(\varepsilon_K\Delta_K = \widetilde O(\varepsilon_K^2).\)
Since \(\varepsilon_K>0\),
\(\Delta_K = f(x_K)-f(x^*) = \widetilde O(\varepsilon_K).\)
\end{proof}

\begin{corollary}
\label[corollary]{cor:det-classical-anchor-rate}
If
\(\varepsilon_n=n^{-1/2}\) and \(\alpha_n=\frac{1}{n+2}\)
for all sufficiently large \(n\), then
\[
        f(x_K)-f(x^*)
        =
        O\!\left(\frac{\log K}{\sqrt K}\right)
        =
        \widetilde O(K^{-1/2}).
\]
\end{corollary}

\begin{proof}
For all sufficiently large \(K\), \(\varepsilon_K=K^{-1/2}\). The finite
prefix contributes only a constant, and hence
\(\sum_{k=1}^{K}\alpha_k=O(\log K)\) and
\(\sum_{k=1}^{K}\varepsilon_k^2=O(\log K)\).
Therefore \eqref{eq:rate-global-det} holds, in fact:
\[
        \frac{
            1+\sum_{k=1}^{K}(\alpha_k+\varepsilon_k^2)
        }{K}
        =
        O\!\left(\frac{\log K}{K}\right)
        =
        \widetilde O(\varepsilon_K^2).
\]

For \eqref{eq:rate-suffix-det}, after enlarging \(C\) to absorb the finite
prefix, \(\alpha_k+\varepsilon_k^2\le C/k\) for every \(k\ge1\). Thus
\cref{lem:polylog-sums}, specifically
\eqref{eq:polylog-suffix-double-sum} with \(s=0\), gives
\[
        \sum_{\ell=1}^{K-1}
        \frac1{\ell(\ell+1)}
        \sum_{k=K-\ell}^{K}
        (\alpha_k+\varepsilon_k^2)
        =
        O\!\left(\frac{\log K}{K}\right)
        =
        \widetilde O(\varepsilon_K^2).
\]
Thus \cref{ass:det-suffix-rate-regime} holds. Inserting the preceding estimates into
the proof of \cref{thm:det-last-iterate-general} gives
\[
        \varepsilon_K\bigl(f(x_K)-f(x^*)\bigr)
        =
        O\!\left(\frac{\log K}{K}\right).
\]
Since \(\varepsilon_K=K^{-1/2}\), we conclude that
\(
        f(x_K)-f(x^*)
        =
        O\!\left(\frac{\log K}{\sqrt K}\right)
        =
        \widetilde O(K^{-1/2}).
\)
\end{proof}

We finally relate the deterministic estimate to the residual bounds known for
classical Halpern fixed-point iterations.
\begin{remark}
If \(0<\eta<2/L\) is fixed and
\(T:=I-\eta\nabla f,\)
then \(T\) is nonexpansive and \(\Fix(T)=S\). The classical Halpern iteration
\(x_{n+1}=\alpha_nu+(1-\alpha_n)Tx_n\)
coincides with the Halpern--gradient method with constant step \(\eta\). In
this fixed-mapping setting, quantitative results are usually stated for the
fixed-point residual \(\norm{x_k-Tx_k}\). For the classical Halpern iteration
initialized at its anchor, that is, in the special case \(u=x_0\), and with
coefficients \(\alpha_k=1/(k+2)\),
Lieder~\cite[Theorem~2.1]{lieder2021} proves, for
\(x^\dagger\in\Fix(T)\), the sharp estimate
\(\frac12\norm{x_k-Tx_k}\le\frac{\norm{x_0-x^\dagger}}{k+1}.\)
Sabach--Shtern~\cite[Lemma~5]{sabach2017first} prove, within their sequential
averaging framework, \(O(1/k)\) estimates for successive differences such as
\(\norm{x_k-x_{k-1}}\) and \(\norm{y_k-x_{k-1}}\), in their notation; in
their bilevel specialization this leads to an \(O(1/k)\) rate for the relevant
inner objective value.

For \(T_\varepsilon:=I-\varepsilon\nabla f\), the residual is exactly
\(\norm{x_k-T_\varepsilon x_k} = \varepsilon\norm{\nabla f(x_k)}.\)
Thus a fixed-point residual bound of order \(O(1/k)\) yields the gradient
residual estimate \(\norm{\nabla f(x_k)}=O(1/(\varepsilon k))\). Along a
bounded trajectory, convexity also gives
\(f(x_k)-f(x^*) \le \inner{\nabla f(x_k)}{x_k-x^*} \le R\norm{\nabla f(x_k)}\)
for some \(R>0\). Hence, in the fixed-step-size regime, the residual estimate
implies the functional last-iterate scale
\[
        f(x_k)-f(x^*)=O\!\left(\frac{1}{\varepsilon k}\right).
\]
For the varying maps
\(T_{\varepsilon_n}:=I-\varepsilon_n\nabla f\), the identity
\(\norm{x_n-T_{\varepsilon_n}x_n}
=\varepsilon_n\norm{\nabla f(x_n)}\) suggests
\(1/(\varepsilon_n n)\) as a heuristic functional scale. The fixed-map
residual estimate alone does not establish this scale, because the operator
now changes with \(n\). Independently, under the suffix assumptions used
above, \cref{cor:det-classical-anchor-rate} attains that benchmark up to a
logarithmic factor: when \(\varepsilon_n=n^{-1/2}\) and
\(\alpha_n=1/(n+2)\) for all sufficiently large \(n\), it gives
\[
        f(x_K)-f(x^*)
        =
        O\!\left(\frac{\log K}{\sqrt K}\right)
        =
        \widetilde O(K^{-1/2}),
\]
that is, the heuristic \(1/(\varepsilon_KK)\) scale, up to logarithmic factors,
while allowing the gradient step to vanish.
\end{remark}

\section{Stochastic convergence}
\label{sec:stoch-convergence}

We use the stochastic setting and noise decomposition of \cref{sec:prelim}.
The bias \(b_n\) may be random and \(\F_n\)-measurable. The quantity
\(\beta_n\) enters the last-iterate estimates, whereas \(L^2\)-boundedness and
pathwise comparison use the scaled squared bias displayed below.

We first record a convenient \(L^2\)-boundedness criterion. It will be used
below in the last-iterate analysis.

\begin{lemma}
\label[lemma]{lem:L2-boundedness-stoch}
Assume the coefficient conditions
\(0<\alpha_n<1, \qquad 0<\varepsilon_n<\frac2L\)
eventually. Suppose that one of the following sufficient conditions holds:
\begin{equation}
\label{eq:L2-summable-sufficient}
        \sum_{n=0}^{\infty}
        \left(
            \E v_n
            +
            \frac{\varepsilon_n^2}{\alpha_n}
            \E\norm{b_n}^2
        \right)
        <+\infty,
\end{equation}
or
\begin{equation}
\label{eq:L2-alpha-sufficient}
        \E v_n
        +
        \frac{\varepsilon_n^2}{\alpha_n}
        \E\norm{b_n}^2
        \le C\alpha_n
        \qquad\text{eventually}.
\end{equation}
Then
\(\sup_{n\ge0}\E\norm{X_n}^2<+\infty.\)
In particular, for every \(p\in S\),
\(\sup_{n\ge0}\E\norm{X_n-p}^2<+\infty,\)
and
\(\sup_{n\ge0}\E\norm{\nabla f(X_n)}^2<+\infty.\)
\end{lemma}

\begin{proof}
Fix \(p\in S\) and work on a tail where the coefficient restrictions hold.
Since \(T_{\varepsilon_n}p=p\), conditional centering, convexity of the squared
norm, nonexpansiveness, and Young's inequality with
\(\delta_n=\alpha_n/[2(1-\alpha_n)]\) give
\[
\begin{aligned}
 \E_n\norm{X_{n+1}-p}^2
 &\le \alpha_n\norm{u-p}^2
 +(1-\alpha_n)
 \bigl(\norm{T_{\varepsilon_n}X_n-p-\varepsilon_nb_n}^2+v_n\bigr)\\
 &\le \alpha_n\norm{u-p}^2
 +\left(1-\frac{\alpha_n}{2}\right)\norm{X_n-p}^2
 +C\frac{\varepsilon_n^2}{\alpha_n}\norm{b_n}^2+v_n .
\end{aligned}
\]
Taking expectations yields
\begin{equation}
\label{eq:L2-boundedness-recursion}
        \E\norm{X_{n+1}-p}^2
        \le
        \left(1-\frac{\alpha_n}{2}\right)
        \E\norm{X_n-p}^2
        +\alpha_n\norm{u-p}^2+C\E v_n
        +
        C\frac{\varepsilon_n^2}{\alpha_n}
        \E\norm{b_n}^2 .
\end{equation}
Set
\(a_n:=\E\norm{X_n-p}^2\), \(R:=2\norm{u-p}^2\), and
\(d_n:=C\E v_n+C\frac{\varepsilon_n^2}{\alpha_n}
\E\norm{b_n}^2.\)
Under \eqref{eq:L2-summable-sufficient},
\(a_{n+1}\le\max\{a_n,R\}+d_n\), so summability of \(d_n\) bounds \(a_n\).
Under \eqref{eq:L2-alpha-sufficient}, eventually \(d_n\le C\alpha_n\), whence
\[
 a_{n+1}\le
 \left(1-\frac{\alpha_n}{2}\right)a_n
 +\frac{\alpha_n}{2}(R+2C)
 \le\max\{a_n,R+2C\}.
\]
The finite prefix is square-integrable by the standing assumptions. Thus
\(\sup_n a_n<+\infty\), and the other claims follow from
\[
        \norm{X_n}^2\le2\norm{X_n-p}^2+2\norm{p}^2,
        \qquad
        \norm{\nabla f(X_n)}\le L\norm{X_n-p}.
\]
\end{proof}

We now prove almost sure convergence of the whole stochastic sequence by a
pathwise comparison with the deterministic Halpern--gradient orbit.

For each \(\omega\in\Omega\), define the
deterministic comparison orbit
\(Y_0(\omega):=X_0(\omega),\)
and
\begin{equation}
\label{eq:det-comparison-orbit}
        Y_{n+1}(\omega)
        =
        \alpha_nu
        +
        (1-\alpha_n)
        \bigl(
            Y_n(\omega)-\varepsilon_n\nabla f(Y_n(\omega))
        \bigr).
\end{equation}
By \cref{thm:det-strong-convergence},
\(Y_n(\omega)\to x^*\) for every \(\omega.\)
The limit is independent of \(Y_0(\omega)\), since the deterministic scheme
selects \(P_S(u)\).

Set
\(e_n:=X_n-Y_n.\)
Subtracting \eqref{eq:det-comparison-orbit} from
\eqref{eq:intro-stoch}, and using the decomposition
\(U_n=b_n+\xi_n\), gives
\[
\begin{aligned}
        e_{n+1}
        &=
        (1-\alpha_n)
        \Bigl[
            e_n
            -
            \varepsilon_n(\nabla f(X_n)-\nabla f(Y_n))
            -
            \varepsilon_n b_n
            -
            \varepsilon_n\xi_n
        \Bigr].
\end{aligned}
\]

We use the following localized Robbins--Siegmund theorem
\cite{robbins1971convergence}.

\begin{lemma}
\label[lemma]{lem:robbins-siegmund}
Let \((\F_n)_{n\ge0}\) be a filtration. Let
\((Z_n)_{n\ge0}\), \((\Theta_n)_{n\ge0}\), and \((\eta_n)_{n\ge0}\)
be nonnegative adapted real-valued processes, with \(Z_n,\Theta_n,\eta_n\in L^1\), such that
\[
        \E[Z_{n+1}\mid\F_n]
        \le
        Z_n-\Theta_n+\eta_n
        \qquad\text{a.s. for all } n\ge0.
\]
Assume that
\(\sum_{n=0}^{\infty}\eta_n<+\infty\) a.s.
Then \(Z_n\) converges almost surely to a finite random variable and
\(\sum_{n=0}^{\infty}\Theta_n<+\infty\)
a.s. In particular, if \((\lambda_n)\subset[0,1]\) and
\(\E[Z_{n+1}\mid\F_n] \le (1-\lambda_n)Z_n+\eta_n \) a.s.,
then \(Z_n\) converges almost surely and
\(\sum_{n=0}^{\infty}\lambda_n Z_n<+\infty\) a.s.
\end{lemma}

\begin{lemma}
\label[lemma]{lem:comparison-estimate}
Assume that
\[
        0<\alpha_n<1,
        \qquad
        \sum_{n=0}^{\infty}\alpha_n=+\infty,
        \qquad
        0<\varepsilon_n<\frac2L
        \quad\text{eventually}.
\]
Assume moreover the comparison summability condition
\begin{equation}
\label{eq:comparison-summability}
        \sum_{n=0}^{\infty}
        \left(
            v_n
            +
            \frac{\varepsilon_n^2}{\alpha_n}\norm{b_n}^2
        \right)
        <+\infty
        \qquad\text{a.s.}
\end{equation}
Then \(e_n\to0\) a.s.
\end{lemma}

\begin{proof}
Choose \(N\) so that the coefficient restrictions hold for \(n\ge N\), and
apply the argument below to this tail. Taking conditional expectations,
\[
\begin{aligned}
        \E_n\norm{e_{n+1}}^2
        &=
        (1-\alpha_n)^2
        \norm{
            T_{\varepsilon_n}X_n
            -
            T_{\varepsilon_n}Y_n
            -
            \varepsilon_nb_n
        }^2                                                     \\
        &\quad
        +
        (1-\alpha_n)^2v_n.
\end{aligned}
\]
We use Young's inequality with
\(\delta_n:=\frac{\alpha_n}{1-\alpha_n}.\)
Since \(T_{\varepsilon_n}\) is nonexpansive,
\[
\begin{aligned}
        \norm{
            T_{\varepsilon_n}X_n
            -
            T_{\varepsilon_n}Y_n
            -
            \varepsilon_nb_n
        }^2
        &\le
        (1+\delta_n)
        \norm{
            T_{\varepsilon_n}X_n
            -
            T_{\varepsilon_n}Y_n
        }^2                                                     \\
        &\quad
        +
        \left(1+\frac1{\delta_n}\right)
        \varepsilon_n^2\norm{b_n}^2                             \\
        &\le
        (1+\delta_n)\norm{e_n}^2
        +
        C\frac{\varepsilon_n^2}{\alpha_n}\norm{b_n}^2 .
\end{aligned}
\]
Multiplying by \((1-\alpha_n)^2\), and using
\((1-\alpha_n)^2(1+\delta_n)=1-\alpha_n,\)
we obtain
\[
        \E_n\norm{e_{n+1}}^2
        \le
        (1-\alpha_n)\norm{e_n}^2
        +
        C\frac{\varepsilon_n^2}{\alpha_n}\norm{b_n}^2
        +
        v_n.
\]

Set
\(Z_n:=\norm{e_n}^2\) and
\(\eta_n:=C\frac{\varepsilon_n^2}{\alpha_n}\norm{b_n}^2+v_n.\)
Then \(Z_n\ge0\), \(\eta_n\ge0\), \(\eta_n\) is
\(\F_n\)-measurable, and
\(\E[Z_{n+1}\mid \F_n]\le (1-\alpha_n) Z_n+\eta_n.\)
By \eqref{eq:comparison-summability},
\(\sum_{n=0}^{\infty}\eta_n<+\infty \) a.s.
The standing \(L^2\) assumptions ensure the required integrability. Applying
\cref{lem:robbins-siegmund} with \(\Theta_n=\alpha_nZ_n\) on the shifted
filtration, \(Z_n\) converges a.s. and
\(\sum_{n=N}^{\infty}\alpha_nZ_n<+\infty\) a.s.
Since \(\sum_{n=N}^{\infty}\alpha_n=+\infty\), the limit of \(Z_n\) must be
zero; hence \(e_n\to0\) a.s.
\end{proof}

\begin{theorem}
\label{thm:stoch-as-convergence}
Assume
the coefficient conditions
\eqref{eq:alpha-base}--\eqref{eq:coupling} and
the comparison summability condition
\eqref{eq:comparison-summability}. Then
\(X_n\to x^*=P_S(u) \)  a.s.
\end{theorem}

\begin{proof}
Let \((Y_n)\) be the deterministic comparison orbit
\eqref{eq:det-comparison-orbit}. By \cref{thm:det-strong-convergence},
\(Y_n\to x^*\) pointwise, while \cref{lem:comparison-estimate} gives
\(X_n-Y_n\to0\) a.s.; hence \(X_n\to x^*\) a.s.
\end{proof}


\section{Stochastic last-iterate estimates}
\label{sec:stoch-last-iterate}
The comparison summability condition yields almost sure convergence but not a
last-iterate rate: the spike examples in
\cref{ex:rare-large-spikes,ex:WS-no-last-iterate-rate} show precisely
that weighted summability may hold and yet a prescribed last-iterate rate may
fail. Conversely, \cref{ex:suffix-no-as-convergence} shows that the
expectation-level suffix condition below does not imply almost sure
convergence. Set \(\Delta_k:=\E[f(X_k)-f(x^*)]\). The bias \(b_n\) remains
\(\F_n\)-measurable and possibly random; if \(b_n=0\), then
\(\beta_n=0\) and \(\xi_n=U_n\).

\begin{assumption}
\label[assumption]{ass:stoch-suffix-regularity}
The sequence \((\varepsilon_n)\) is eventually nonincreasing. With
\(\theta_n := \alpha_n+\varepsilon_n^2+\E v_n+\beta_n,\)
one has
\begin{align}
        \frac{
            1+\sum_{k=1}^{K}\theta_k
        }{K}
        &=
        \widetilde O(\varepsilon_K^2),
        \label{eq:stoch-rate-global}
        \\
        \sum_{\ell=1}^{K-1}
        \frac1{\ell(\ell+1)}
        \sum_{k=K-\ell}^{K}\theta_k
        &=
        \widetilde O(\varepsilon_K^2).
        \label{eq:stoch-rate-suffix}
\end{align}
\end{assumption}

\begin{theorem}
\label{thm:stoch-last-iterate-general}
Assume the coefficient conditions
\eqref{eq:alpha-base}--\eqref{eq:coupling}, the hypotheses of
\cref{lem:L2-boundedness-stoch}, and the stochastic suffix-regularity condition
\cref{ass:stoch-suffix-regularity}. Then
\(\E[f(X_K)-f(x^*)] = \widetilde O(\varepsilon_K).\)
\end{theorem}

\begin{proof}
Fix an \(\F_k\)-measurable \(z\) with bounded second moment. Conditional
centering of \(\xi_k\), convexity of the squared norm and of \(f\), and
expansion of the square give the following estimate. For \(z=x^*\), and
uniformly for \(z=X_m\) when \(k\ge m\), the \(L^2\)-bounds of
\cref{lem:L2-boundedness-stoch} and Cauchy--Schwarz control the bias by
\[
 \varepsilon_k\left|
 \E\inner{X_k-z-\varepsilon_k\nabla f(X_k)}{b_k}\right|
 +\varepsilon_k^2\E\norm{b_k}^2
 \le C\beta_k.
\]
Consequently,
\begin{equation}
\label{eq:stoch-one-step-comparator}
\begin{aligned}
        &2(1-\alpha_k)\varepsilon_k
        \E\bigl[f(X_k)-f(z)\bigr]
        +
        \E\norm{X_{k+1}-z}^2                              \\
        &\le
        (1-\alpha_k)\E\norm{X_k-z}^2
        +
        \alpha_k\E\norm{u-z}^2
        +
        C\varepsilon_k^2
        +
        C\E v_k
        +
        C\beta_k .
\end{aligned}
\end{equation}

Set \(\Delta_k:=\E[f(X_k)-f(x^*)]\). Smoothness and
\cref{lem:L2-boundedness-stoch} give \(0\le\Delta_k\le C\). Choose \(N_0\) so
that, for \(k\ge N_0\), all eventual restrictions hold,
\(\alpha_k\le1/2\), \(\varepsilon_k\le1\), and \((\varepsilon_k)\) is
nonincreasing. Summing \eqref{eq:stoch-one-step-comparator} with \(z=x^*\)
and absorbing the finite prefix gives
\begin{equation}
\label{eq:average-stoch-rate}
        \frac1K
        \sum_{k=1}^{K}
        \varepsilon_k\Delta_k
        \le
        C
        \frac{
            1+\sum_{k=1}^{K}\theta_k
        }{K}.
\end{equation}

Fix \(m=K-\ell\ge N_0\) and set
\(D_{k,m}:=\E[f(X_k)-f(X_m)]=\Delta_k-\Delta_m\), so
\(|D_{k,m}|\le C\). In \eqref{eq:stoch-one-step-comparator} take \(z=X_m\)
and sum over \(k=m,\ldots,K\). The distance terms telescope to a nonpositive
quantity and the anchor terms are bounded by \(C\sum_{k=m}^K\alpha_k\).
After dividing by \(2\) and adding the signed remainder,
\begin{equation}
\label{eq:suffix-stoch-estimate}
\begin{aligned}
 \sum_{k=m}^{K}\varepsilon_kD_{k,m}
 &=
 \sum_{k=m}^{K}(1-\alpha_k)\varepsilon_kD_{k,m}
 +\sum_{k=m}^{K}\alpha_k\varepsilon_kD_{k,m}\\
 &\le C\sum_{k=m}^{K}\theta_k
 +C\sum_{k=m}^{K}\alpha_k\varepsilon_k
 \le C\sum_{k=m}^{K}\theta_k .
\end{aligned}
\end{equation}
Since \(\varepsilon_k\le\varepsilon_m\) and \(\Delta_m\ge0\), this implies
\[
 \sum_{k=m}^{K}(\varepsilon_k\Delta_k-\varepsilon_m\Delta_m)
 \le C\sum_{k=m}^{K}\theta_k.
\]
For the finitely many \(m<N_0\), nonnegativity and
\eqref{eq:average-stoch-rate} bound each inner suffix sum by
\(C(1+\sum_{k=1}^K\theta_k)\). Their suffix weights are \(O(K^{-2})\).
Thus the identity \eqref{eq:suffix-averaging-identity} in
\cref{lem:suffix-averaging} yields
\[
\varepsilon_K\Delta_K
\le
\frac1K\sum_{k=1}^{K}\varepsilon_k\Delta_k
+C\sum_{\ell=1}^{K-1}\frac1{\ell(\ell+1)}
  \sum_{k=K-\ell}^{K}\theta_k
+\frac{C}{K^2}\left(1+\sum_{k=1}^{K}\theta_k\right)
=\widetilde O(\varepsilon_K^2),
\]
by \eqref{eq:stoch-rate-global}--\eqref{eq:stoch-rate-suffix}. Dividing by
\(\varepsilon_K\) completes the proof.
\end{proof}

\begin{corollary}
\label[corollary]{cor:stoch-classical-anchor-rate}
Assume the unbiased case \(b_n=0\), and suppose that
\(\E[\|U_n\|^2\mid\F_n]\le\sigma^2 \) a.s.
Assume moreover that, for all sufficiently large \(n\),
\(\varepsilon_n=n^{-1/2}\) and \(\alpha_n=\frac1{n+2}\).
Then
\[
        \E[f(X_K)-f(x^*)]
        =
        O\left(\frac{\log K}{\sqrt K}\right)
        =
        \widetilde O(K^{-1/2}).
\]
\end{corollary}

\begin{proof}
Here \(\beta_n=0\), \(\E v_n\le\sigma^2/n\), and hence
\(\theta_n\le C/n\) and \(\E v_n\le C\alpha_n\) eventually. Thus
\eqref{eq:L2-alpha-sufficient} holds, while
\cref{lem:polylog-sums}, specifically
\eqref{eq:polylog-suffix-double-sum} with \(s=0\), verifies both suffix
conditions with the quantitative bound \(O((\log K)/K)\). Substitution in the
proof of \cref{thm:stoch-last-iterate-general} gives
\(\varepsilon_K\E[f(X_K)-f(x^*)]=O((\log K)/K)\), and
\(\varepsilon_K=K^{-1/2}\) completes the proof.
\end{proof}

\begin{remark}
\label{rem:bv-rate-not-as}
For \(\varepsilon_n=n^{-1/2}\) and
\(\alpha_n=(\log n)/n\), bounded variance gives \(\E v_n=O(1/n)\), which is
enough for the last-iterate bound but not for
\(\sum_n v_n<+\infty\) a.s. The distinction is realized by
\cref{ex:suffix-no-as-convergence}.
\end{remark}

\section{Bounded higher moments as an anti-spike condition}
\label{sec:higher-moments}

In the unbiased case, pathwise comparison requires
\(\sum_n v_n<+\infty\) a.s., whereas the last-iterate theorem only controls
expectation-level suffix sums. At the critical scale
\(\varepsilon_n=n^{-1/2}\), \(\alpha_n=(\log n)/n\), bounded conditional
second moments give \(v_n=O(1/n)\): the rate applies, but comparison need not;
see \cref{ex:suffix-no-as-convergence}. The block argument below instead uses
a uniform conditional \(p\)-moment, \(p>4\), to control martingale
perturbations by Burkholder and Borel--Cantelli. The restriction \(p>4\) is
used only for almost sure convergence; the expected rate uses second moments
alone, and no optimality of this threshold is claimed.

\begin{lemma}
\label[lemma]{lem:local-nonexpansive-block-estimate}
Let \(p>2\). Let \(a<b\) be integers, and let \((z_k)_{k=a}^{b}\) be an
adapted \(\HH\)-valued process satisfying
\(z_a=0\)
and
\(z_{k+1} = q_k \bigl( R_k z_k-\varepsilon_k U_k \bigr), \qquad k=a,\ldots,b-1,\)
where \((q_k)_{k=a}^{b-1}\) and \((\varepsilon_k)_{k=a}^{b-1}\) are
deterministic sequences satisfying \(0\le q_k\le1\) and
\(\varepsilon_k\ge0\), and \(R_k:\HH\to\HH\) is a possibly random map such
that, for every \(x\in\HH\), \(R_kx\) is
\(\F_k\)-measurable, \(R_k\) is nonexpansive, \(R_k0=0\), and
\(R_kz_k\) is \(\F_k\)-measurable. Assume
that \(U_k\) is \(\F_{k+1}\)-measurable and
\(\E[U_k\mid\F_k]=0\), \(k=a,\ldots,b-1\),
and that, for some constant \(C_p>0\),
\(\E[\|U_k\|^p\mid\F_k]\le C_p\) a.s.
Then there exists a constant \(C\), depending only on \(p\) and \(C_p\), such
that
\[
        \E
        \sup_{a\le n\le b}\|z_n\|^p
        \le
        C
        \left(
            \sum_{k=a}^{b-1}\varepsilon_k^2
        \right)^{p/2}.
\]
\end{lemma}

\begin{proof}
Set
\(H:=\sup_{a\le n\le b}\|z_n\|^2\).
For \(N>0\), stop at
\(
        \tau_N:=\inf\{n\in\{a,\ldots,b\}:\|z_n\|\ge N\},
\)
with \(\tau_N=b+1\) if the set is empty, and write
\(\mathbf 1_k^N:=\mathbf 1_{\{k<\tau_N\}}\). Then
\(\mathbf 1_k^N\) is \(\F_k\)-measurable.

Define
\(H_N := \sup_{a\le j\le b}\|z_{j\wedge\tau_N}\|^2\)
and
\(\Gamma_N := \sum_{k=a}^{b-1} \mathbf 1_k^N q_k^2\varepsilon_k^2\|U_k\|^2\).
Also define
\[
        M_j^N
        :=
        -2
        \sum_{k=a}^{j-1}
        \mathbf 1_k^N q_k^2\varepsilon_k
        \inner{R_kz_k}{U_k},
        \qquad j=a,\ldots,b.
\]
Then \((M_j^N)_{j=a}^{b}\) is a real-valued martingale. Indeed, the increment is
\(\Delta M_{k+1}^N = -2\mathbf 1_k^Nq_k^2\varepsilon_k \inner{R_kz_k}{U_k},\)
where \(\mathbf 1_k^Nq_k^2\varepsilon_k R_kz_k\) is
\(\F_k\)-measurable. Hence
\[
        \E[\Delta M_{k+1}^N\mid\F_k]
        =
        -2\mathbf 1_k^Nq_k^2\varepsilon_k
        \inner{R_kz_k}{\E[U_k\mid\F_k]}
        =
        0.
\]
Moreover, on \(\{k<\tau_N\}\),
\(\|R_kz_k\|\le \|z_k\|<N,\)
so the stopped increments are integrable.

Since \(R_k0=0\) and \(R_k\) is nonexpansive,
\(\|R_kz_k\| \le \|z_k\|.\)
Therefore,
\[
\begin{aligned}
        \|z_{k+1}\|^2
        &=
        q_k^2
        \|R_kz_k-\varepsilon_kU_k\|^2                            \\
        &\le
        \|z_k\|^2
        -
        2q_k^2\varepsilon_k
        \inner{R_kz_k}{U_k}
        +
        q_k^2\varepsilon_k^2\|U_k\|^2 .
\end{aligned}
\]
Summing this inequality up to the stopped time gives, for
\(j=a,\ldots,b\),
\[
        \|z_{j\wedge\tau_N}\|^2 \le \Gamma_N+M_j^N.
\]
Hence
\[
        H_N \le \Gamma_N+\sup_{a\le j\le b}|M_j^N|.
\]
Let \(r:=\frac p2>1\). Then
\[
        \E H_N^r
        \le C_r\E\Gamma_N^r
        +C_r\E\sup_{a\le j\le b}|M_j^N|^r.
\]
By the Burkholder--Davis--Gundy inequality for real-valued martingales,
\[
\begin{aligned}
        \E\sup_{a\le j\le b}|M_j^N|^r
        &\le
        C_r
        \E
        \left(
            \sum_{k=a}^{b-1}
            (\Delta M_{k+1}^N)^2
        \right)^{r/2}
        \le
        C_r
        \E
        \left(
            \sum_{k=a}^{b-1}
            \mathbf 1_k^Nq_k^4\varepsilon_k^2
            \inner{R_kz_k}{U_k}^2
        \right)^{r/2}.
\end{aligned}
\]
Moreover,
\(|\inner{R_kz_k}{U_k}| \le \|R_kz_k\|\,\|U_k\| \le \|z_k\|\,\|U_k\|.\)
On \(\{k<\tau_N\}\), \(\|z_k\|^2\le H_N\). Thus
\[
\begin{aligned}
        \sum_{k=a}^{b-1}
        \mathbf 1_k^Nq_k^4\varepsilon_k^2
        \inner{R_kz_k}{U_k}^2
        &\le
        H_N
        \sum_{k=a}^{b-1}
        \mathbf 1_k^Nq_k^2\varepsilon_k^2\|U_k\|^2                =
        H_N\Gamma_N .
\end{aligned}
\]
Consequently,
\(\E\sup_{a\le j\le b}|M_j^N|^r \le C_r \E(H_N\Gamma_N)^{r/2}.\)
By Young's inequality, for every \(\eta>0\),
\((H_N\Gamma_N)^{r/2} = H_N^{r/2}\Gamma_N^{r/2} \le \eta H_N^r+C_{\eta,r}\Gamma_N^r.\)
Choosing \(\eta>0\) sufficiently small and absorbing the corresponding term on
the left-hand side, we obtain
\(\E H_N^r \le C_r\E\Gamma_N^r.\)
Since \(0\le \Gamma_N\le \Gamma\), where
\(\Gamma := \sum_{k=a}^{b-1} q_k^2\varepsilon_k^2\|U_k\|^2,\)
we have
\(\E H_N^r \le C_r\E\Gamma^r.\)
Finally, since \(q_k\le1\), Minkowski's inequality in \(L^r\) gives
\[
        \|\Gamma\|_{L^r}
        \le
        \sum_{k=a}^{b-1}\varepsilon_k^2\bigl\|\|U_k\|^2\bigr\|_{L^r}
        =
        \sum_{k=a}^{b-1}\varepsilon_k^2\bigl(\E\|U_k\|^p\bigr)^{2/p}
        \le
        C\sum_{k=a}^{b-1}\varepsilon_k^2.
\]
Therefore
\(
        \E\Gamma^r
        \le
        C
        \left(
            \sum_{k=a}^{b-1}\varepsilon_k^2
        \right)^r.
\)
The bound is uniform in \(N\). Since the interval \(\{a,\ldots,b\}\) is finite,
we have \(H_N\to H\) almost surely as \(N\to\infty\). Therefore, by Fatou's
lemma,
\[
        \E H^r
        \le
        \liminf_{N\to\infty}\E H_N^r
        \le
        C
        \left(
            \sum_{k=a}^{b-1}\varepsilon_k^2
        \right)^r.
\]
Since \(H^r=\sup_{a\le n\le b}\|z_n\|^p\), the proof is complete.
\end{proof}

\begin{theorem}
\label{thm:moment-as-and-rate}
Assume the unbiased case
\(b_n=0, \qquad U_n=\xi_n, \qquad \E[U_n\mid\F_n]=0.\)
Assume that, for some \(p>4\), there exists \(C_p>0\) such that
\begin{equation}
\label{eq:moment-p-greater-four-main}
        \E[\|U_n\|^{p}\mid\F_n]
        \le C_p
        \qquad\text{a.s. for all }n.
\end{equation}
Assume moreover that, for all sufficiently large \(n\),
\(\varepsilon_n=n^{-1/2},\) \(\alpha_n=\frac{\log n}{n}.\)
Then
\(X_n\to x^*=P_S(u)\)  a.s.
Furthermore,
\[
        \E[f(X_K)-f(x^*)]
        =
        O\left(
            \frac{(\log K)^2}{\sqrt K}
        \right).
\]
\end{theorem}
\begin{proof}
Constants denoted by \(C\) may change from line to line. Conditional Jensen's
inequality and \eqref{eq:moment-p-greater-four-main} give
\(\E[\|U_n\|^2\mid\F_n]\le C\), hence \(v_n\le C/n\) eventually.
Let \((Y_n)\) be the deterministic comparison orbit
\eqref{eq:det-comparison-orbit}, with \(Y_0=X_0\). By
\cref{thm:det-strong-convergence}, \(Y_n\to x^*\) a.s. Set
\(e_n:=X_n-Y_n\) and
\(R_nz := T_{\varepsilon_n}(Y_n+z)-T_{\varepsilon_n}(Y_n).\)
Eventually \(R_n\) is nonexpansive, \(R_n0=0\), and
\begin{equation}
\label{eq:error-recursion-short-moment-proof}
        e_{n+1}
        =
        (1-\alpha_n)
        \bigl(
            R_ne_n-\varepsilon_nU_n
        \bigr).
\end{equation}

Choose \(\nu_0\) so that the preceding facts hold for \(n\ge\nu_0\), and define
\[
        \nu_{j+1}
        :=
        \min
        \left\{
            n>\nu_j:
            \sum_{k=\nu_j}^{n-1}\alpha_k\ge1
        \right\}, \qquad I_j:=\{\nu_j,\ldots,\nu_{j+1}-1\}.
\]
The endpoints exist because \(\sum_n\alpha_n=\infty\), and minimality gives
\(1\le\sum_{k\in I_j}\alpha_k<2\).

On each block \(I_j\), define the zero-started perturbation process
\(w_{\nu_j}^{(j)}=0,\)
\[
        w_{k+1}^{(j)}
        =
        (1-\alpha_k)
        \bigl(
            R_kw_k^{(j)}-\varepsilon_kU_k
        \bigr),
        \qquad k=\nu_j,\ldots,\nu_{j+1}-1.
\]
The deterministic coefficients, adaptedness, and continuity of
\(T_{\varepsilon_k}\) verify the measurability hypotheses of
\cref{lem:local-nonexpansive-block-estimate}.
Let
\(Z_j := \sup_{\nu_j\le n\le\nu_{j+1}} \|w_n^{(j)}\|.\)
By \cref{lem:local-nonexpansive-block-estimate},
\[
        \E Z_j^p
        \le
        C
        \left(
            \sum_{k=\nu_j}^{\nu_{j+1}-1}\varepsilon_k^2
        \right)^{p/2}.
\]
Since \(\varepsilon_k^2=1/k=\alpha_k/\log k\),
\[
        \sum_{k=\nu_j}^{\nu_{j+1}-1}\varepsilon_k^2
        \le
        \frac1{\log \nu_j}
        \sum_{k=\nu_j}^{\nu_{j+1}-1}\alpha_k
        \le
        \frac{2}{\log \nu_j}.
\]
Also,
\[
        j
        \le
        \sum_{k=\nu_0}^{\nu_j-1}\alpha_k
        =
        \sum_{k=\nu_0}^{\nu_j-1}\frac{\log k}{k}
        \le
        C+\frac12(\log \nu_j)^2.
\]
Thus \(\log\nu_j\ge c\sqrt j\) and
\(\E Z_j^p\le Cj^{-p/4}\). Since \(p>4\), for every \(r\in\mathbb N\),
\[
 \sum_j\PP(Z_j>1/r)\le r^p\sum_j\E Z_j^p<\infty.
\]
Borel--Cantelli, applied to the countable family \(r\in\mathbb N\), gives
\(Z_j\to0\) a.s.

On \(I_j\), nonexpansiveness gives
\[
        \|e_n-w_n^{(j)}\|
        \le
        \left(
            \prod_{i=\nu_j}^{n-1}(1-\alpha_i)
        \right)
        \|e_{\nu_j}\|.
\]
At the endpoint, the product is at most
\(\exp(-\sum_{i\in I_j}\alpha_i)\le e^{-1}\); hence
\(\|e_{\nu_{j+1}}\|\le e^{-1}\|e_{\nu_j}\|+Z_j\), so
\(e_{\nu_j}\to0\) a.s. Moreover,
\(\sup_{n\in I_j}\|e_n\|\le\|e_{\nu_j}\|+Z_j\), proving
\(e_n\to0\) and therefore \(X_n\to x^*\) a.s.

Finally, \(\beta_n=0\), \(\E v_n\le C/n\), and
\(\theta_n\le C(\log n)/n\). Thus \eqref{eq:L2-alpha-sufficient} holds, and
\cref{lem:polylog-sums} with \(s=1\), inserted quantitatively into
\cref{thm:stoch-last-iterate-general}, gives the asserted expected rate.
\end{proof}


\section{Applications and motivation of the assumptions}
\label{sec:appl}

The update scales the centered oracle error by \(\varepsilon_n\); hence its
relevant variance is the effective displacement \(v_n\) in
\eqref{eq:def-vn}, not the raw conditional variance alone. This distinction is
standard in stochastic approximation~\cite{kushner2003stochastic} and is
essential for zeroth-order estimators, whose observation noise is amplified as
the differencing radius vanishes
\cite{kiefer1952stochastic,spall1992multivariate,nesterov2017random}.

\subsection{Finite-difference framework}
\label{subsec:finite-difference-application}

Let \(\HH=\mathbb R^d\). The query points \(X_n\pm\rho_ne_i\) are
\(\F_n\)-measurable, while the conditionally independent noises
\(\eta_{n,i}^{\pm}\) are \(\F_{n+1}\)-measurable and satisfy
\(\En\eta_{n,i}^{\pm}=0\) and
\(\En|\eta_{n,i}^{\pm}|^2\le\sigma^2\). With
\(Y_{n,i}^{\pm}=f(X_n\pm\rho_ne_i)+\eta_{n,i}^{\pm}\), consider
\begin{equation}
\label{eq:coordinate-finite-difference-estimator}
        \widehat g_{n,i}
        =
        \frac{Y_{n,i}^+-Y_{n,i}^-}{2\rho_n},
        \qquad i=1,\ldots,d,
\end{equation}
where \(\rho_n\downarrow0\) is the differencing radius; see
\cite{kiefer1952stochastic,berahas2019derivative,more2011estimating}.
Writing \(\widehat g_n=\nabla f(X_n)+b_n+\xi_n\), set
\[
 b_n:=\En\widehat g_n-\nabla f(X_n),
 \qquad \xi_n:=\widehat g_n-\En\widehat g_n.
\]
Then \(b_n\) is \(\F_n\)-measurable and \(\En\xi_n=0\).

Suppose that the oracle satisfies
\begin{equation}
\label{eq:common-zeroth-order-bounds}
        \E[\|\xi_n\|^2\mid\F_n]
        \le C\kappa_n\rho_n^{-2},
        \qquad
        \|b_n\|\le C\rho_n^q
        \quad\text{a.s.},
\end{equation}
where \(q\in\{1,2\}\) and \((\kappa_n)\) is a deterministic positive sequence
encoding its variance scale. Under only a Lipschitz gradient, both one-sided
and symmetric differences generally have \(q=1\); a symmetric difference has
\(q=2\), for example, under a Lipschitz Hessian. The unreplicated coordinate
estimator has \(\kappa_n=1\), with dimension and \(\sigma^2\) absorbed into
\(C\). From \eqref{eq:common-zeroth-order-bounds},
\[
        v_n\le C\varepsilon_n^2\kappa_n\rho_n^{-2},
        \qquad
        \beta_n
        \le C\varepsilon_n\rho_n^q
        +C\varepsilon_n^2\rho_n^{2q}.
\]
Thus the pathwise comparison theorem applies if
\begin{align}
\label{eq:WS-zeroth-order}
        \sum_{n=0}^{\infty}
        \varepsilon_n^2\kappa_n\rho_n^{-2}
        &<+\infty,\\
\label{eq:pathwise-bias-summability-zeroth-order}
        \sum_{n=0}^{\infty}
        \frac{\varepsilon_n^2}{\alpha_n}\rho_n^{2q}
        &<+\infty.
\end{align}
The additional series \(\sum_n\varepsilon_n\rho_n^q<\infty\) implies
\(\sum_n\beta_n<\infty\), but is not required for the almost sure comparison.
Last-iterate rates instead require the global and suffix conditions of
\cref{ass:stoch-suffix-regularity} for
\(\theta_n=\alpha_n+\varepsilon_n^2+\E v_n+\beta_n\).

\begin{corollary}
\label[corollary]{cor:finite-difference-as}
Assume \eqref{eq:alpha-base}--\eqref{eq:coupling}. Let a finite-difference
oracle with radius \(\rho_n\downarrow0\) satisfy
\eqref{eq:common-zeroth-order-bounds} for some \(q\in\{1,2\}\) and some
deterministic positive sequence \((\kappa_n)\). If
\eqref{eq:WS-zeroth-order} and
\eqref{eq:pathwise-bias-summability-zeroth-order} hold, then the corresponding
HalpernSGD iterates satisfy
\(X_n\to P_S(u)\) a.s.
\end{corollary}

\begin{proof}
The two series imply \eqref{eq:comparison-summability}; apply
\cref{thm:stoch-as-convergence}.
\end{proof}

For instance, in the unreplicated symmetric case \(q=2\), \(\kappa_n=1\),
choose any admissible finite prefix and, eventually, set
\[
 \varepsilon_n=(n+1)^{-3/4},\quad
 \alpha_n=(n+1)^{-4/5},\quad \rho_n=(n+1)^{-1/5}.
\]
The admissible raw scale \(\rho_n^{-2}=(n+1)^{2/5}\) diverges, while
\(\varepsilon_n^2\rho_n^{-2}\), \(\varepsilon_n\rho_n^2\), and
\(\alpha_n^{-1}\varepsilon_n^2\rho_n^4\) have orders
\((n+1)^{-11/10}\), \((n+1)^{-23/20}\), and \((n+1)^{-3/2}\), respectively.

\subsection{Random directions}
\label{subsec:random-direction-application}

SPSA~\cite{spall1992multivariate} samples an
\(\F_{n+1}\)-measurable direction \(\Delta_n\), independent of \(\F_n\),
with nonzero coordinates. Given
\(\mathcal G_n=\F_n\vee\sigma(\Delta_n)\), let the fresh observation errors be
centered with bounded variance conditional on \(\mathcal G_n\). For
\(Y_n^\pm=f(X_n\pm\rho_n\Delta_n)+\eta_n^\pm\), set
\[
        \widehat g_{n,i}
        =
        \frac{Y_n^+-Y_n^-}
        {2\rho_n(\Delta_n)_i},
        \qquad i=1,\ldots,d.
\]
This uses two evaluations, but direction sampling generally adds a term of
order \(\|\nabla f(X_n)\|^2\) to the second moment; thus
\eqref{eq:common-zeroth-order-bounds} is not automatic. SPSA and related
two-point estimators
\cite{nesterov2017random,duchi2015optimal,ghadimi2013stochastic,shamir2017optimal}
are applications of the abstract result whenever their bias has the stated
order and, along the generated orbit,
\[
        \E[\|\xi_n\|^2\mid\F_n]
        \le C_d\kappa_n\rho_n^{-2}.
\]
For an unreplicated oracle, \(\kappa_n=1\) is sufficient under appropriate
direction and inverse-direction moments and an a.s. uniform bound on
\(\|\nabla f(X_n)\|\). In particular, SPSA is only a conditional application,
not an oracle model valid for every smooth convex objective.

\subsection{Independent replications}
\label{subsec:noisy-simulation-application}

If each function value is averaged over \(m_n\) conditionally independent
replications, then
\[
        \E[\|\xi_n\|^2\mid\F_n]
        \le \frac{C}{m_n\rho_n^2},
        \qquad
        v_n\le\frac{C\varepsilon_n^2}{m_n\rho_n^2}.
\]
The bias is unchanged, so \(\kappa_n=m_n^{-1}\).

\begin{corollary}
\label[corollary]{cor:replicated-finite-difference-as}
Assume \eqref{eq:alpha-base}--\eqref{eq:coupling}. Let
\(\rho_n\downarrow0\), let \(q\in\{1,2\}\), and let
\((m_n)\subset\mathbb N\) be positive. Suppose that the noise in each
function-value observation is conditionally centered and has conditional
variance at most \(\sigma^2\), that the \(m_n\) replications are conditionally independent
given \(\F_n\), and that the replicated finite-difference estimator
has bias \(\|b_n\|\le C\rho_n^q\) a.s. If
\begin{equation}
\label{eq:replicated-variance-summability}
        \sum_{n=0}^{\infty}
        \frac{\varepsilon_n^2}{m_n\rho_n^2}<+\infty
\end{equation}
and \eqref{eq:pathwise-bias-summability-zeroth-order} holds, then
\(X_n\to P_S(u)\) a.s.
\end{corollary}

\begin{proof}
Apply \cref{cor:finite-difference-as} with \(\kappa_n=m_n^{-1}\).
\end{proof}

\subsection{A stochastic Nesterov-type interpretation}
\label{subsec:stochastic-nesterov-application}

Tran-Dinh~\cite{trandinh2024halpern} related anchored deterministic Halpern
iterations to Nesterov-type recursions. The same algebra applies here. Define
\(Z_{n+1}=X_n-\varepsilon_nG_n\) and, for \(n\ge1\),
\begin{equation}
    \label{eq:nesterov-coefficients}
 \vartheta_n:=\frac{\alpha_n(1-\alpha_{n-1})}{\alpha_{n-1}},
 \qquad \upsilon_n:=\frac{\alpha_n}{\alpha_{n-1}}.
\end{equation}
Eliminating the anchor gives the equivalent recursion
\begin{equation}
    \label{eq:corrected-nesterov-recursion}
\left\{
\begin{aligned}
        Z_{n+1}
        &=
        X_n-\varepsilon_nG_n,
        \\
        X_{n+1}
        &=
        Z_{n+1}
        +
        \vartheta_n(Z_{n+1}-Z_n)
        +
        \upsilon_n(X_n-Z_{n+1}),
\end{aligned}
\right.
\end{equation}
with compatible initialization
\begin{equation}
\label{eq:corrected-nesterov-initialization}
    Z_1=X_0-\varepsilon_0G_0, \qquad X_1=\alpha_0u+(1-\alpha_0)Z_1.
\end{equation}
Here \(\upsilon_n(X_n-Z_{n+1})=\upsilon_n\varepsilon_nG_n\) is the anchoring
correction. Thus this is a particular algebraic representation of HalpernSGD,
not a claim about arbitrary stochastic momentum~\cite{nesterov1983method}.

\begin{corollary}
\label[corollary]{cor:corrected-nesterov-as}
Assume the standing range \(0<\varepsilon_n\le2/L\) and the coefficient
conditions \eqref{eq:alpha-base}--\eqref{eq:coupling}. Let
\((X_n,Z_n)\) be generated by the corrected stochastic Nesterov recursion
\eqref{eq:corrected-nesterov-recursion}, with \(\vartheta_n\) and
\(\upsilon_n\) as in \eqref{eq:nesterov-coefficients} and compatible
initialization as in \eqref{eq:corrected-nesterov-initialization}.
Let \(U_n=b_n+\xi_n\), with \(b_n\) \(\F_n\)-measurable and
\(\En\xi_n=0\), and define
\(D_n=\varepsilon_nb_n\),
\(\Delta M_{n+1}=\varepsilon_n\xi_n\), and
\(s_n=\En\|\Delta M_{n+1}\|^2\). If
\[
        \sum_{n=0}^{\infty}
        \left(
            s_n
            +
            \frac{\|D_n\|^2}{\alpha_n}
        \right)
        <+\infty
        \qquad\text{a.s.},
\]
then
\((X_n,Z_n)\to(P_S(u),P_S(u))\) a.s.\ in \(\HH\times\HH\).
\end{corollary}

\begin{proof}
Since the summability assumption \eqref{eq:comparison-summability} holds,
\cref{thm:stoch-as-convergence} yields
\(X_n\to x^*=P_S(u)\) a.s. Moreover, the same summability assumption
implies that \(D_n\to0\).
For \(A_n^\eta=\{\|\Delta M_{n+1}\|>\eta\}\in\F_{n+1}\),
conditional Markov yields
\(\mathbb P(A_n^\eta\mid\F_n)\le\eta^{-2}s_n\); conditional
Borel--Cantelli therefore gives \(\Delta M_{n+1}\to0\) a.s. Finally,
\(\varepsilon_n\nabla f(X_n)\to0\) and
\(Z_{n+1}=X_n-\varepsilon_n\nabla f(X_n)-D_n-\Delta M_{n+1}\).
\end{proof}

\section{Comparison with related work}
\label{sec:comparison}

The closest results address different combinations of problem class, oracle
structure, and output criterion. Our pathwise assumptions yield strong
selection of \(P_S(u)\), whereas the expectation-level suffix assumptions yield
last-iterate rates; neither regime dominates the alternatives below.

\subsection{Alacaoglu--Malitsky--Wright}
\label{subsec:comparison-AMW}

Alacaoglu--Malitsky--Wright~\cite{alacaoglu2025weakervar} treat constrained,
possibly nonsmooth convex minimization in \(\mathbb R^d\). Their oracle is
unbiased, \(\E[\widetilde g(x)]\in\partial f(x)\), and obeys the
Blum--Gladyshev bound
\[
 \E\|\widetilde g(x)\|^2\le B^2\|x-x_0\|^2+G^2.
\]
Their projected iteration is
\[
 x_{k+1}=P_{\mathcal X}\!\left(
   \beta_kx_0+(1-\beta_k)x_k-\tau_k\widetilde g(x_k)\right),
 \qquad \beta_k=(k+2)^{-1},\quad \tau_k\asymp k^{-1/2},
\]
and achieves an \(O(k^{-1/2})\) weighted-average rate and a horizon-free
\(\widetilde O(k^{-1/2})\) last-iterate rate.

At \(\alpha_n=(n+2)^{-1}\), \(\varepsilon_n=n^{-1/2}\) eventually, our
bounded-variance last-iterate rate is \(O(\log K/\sqrt K)\), matching their
soft order. Their pointwise full-oracle growth bound neither implies nor is
implied by our along-trajectory comparison summability: persistent critical
bounded variance gives \(v_n=O(n^{-1})\), whereas our scaling can accommodate
diverging raw variance. Their setting is broader in nonsmoothness and
constraints; ours allows an \(\F_n\)-measurable bias and yields strong convergence to the
anchor-selected minimizer under its separate pathwise regime.

\subsection{Bravo--Contreras}
\label{subsec:comparison-Bravo-Contreras}

Bravo--Contreras~\cite{bravo2024stochastic} study fixed nonexpansive or
contractive maps in finite-dimensional normed spaces through an unbiased
minibatch oracle. With their convention
\[
 x_n=(1-\beta_n^{\mathrm{BC}})x_0
      +\beta_n^{\mathrm{BC}}(Tx_{n-1}+U_n),
 \qquad \beta_n^{\mathrm{BC}}\to1,
\]
for \(\beta_n^{\mathrm{BC}}=n/(n+1)\) and
\(\sup_n\E\|Tx_n-x_0\|\le\kappa\), the condition
\(\sum_{n\le N}n\E\|U_n\|=\widetilde O(1)\) yields
\(\E\|x_N-Tx_N\|=\widetilde O(N^{-1})\). Under uniformly bounded Euclidean
oracle variance,
their choice \(k_n=n^4\) gives \(\widetilde O(\epsilon^{-5})\) oracle
complexity. Under the usual conditions on \((\beta_n^{\mathrm{BC}})\), if
\(\sum_n\E\|U_n\|<\infty\), the residual converges a.s. and all cluster points
are fixed; whole-sequence a.s. convergence to a random fixed point is obtained
when the normed space is smooth.

Their criterion is the residual of a fixed map. Here
\(T_{\varepsilon_n}=I-\varepsilon_n\nabla f\) varies with \(n\) and its
residual is \(\varepsilon_n\|\nabla f(x)\|\), so their complexity does not
translate into our objective bound without an error condition. Conversely, our
step-scaled variance condition can avoid minibatching and selects \(P_S(u)\),
but does not provide their general fixed-map oracle complexity.

\subsection{Cai--Song--Guzm\'an--Diakonikolas}
\label{subsec:comparison-Cai-Song-Guzman-Diakonikolas}

Cai--Song--Guzm\'an--Diakonikolas~\cite{cai2022stochastic} treat
finite-dimensional cocoercive and Lipschitz-monotone equations. Their oracle is
unbiased, has uniformly bounded variance, and is mean-square Lipschitz; its
multipoint form can evaluate two arguments with the same sample. Exploiting
this structure, their PAGE--Halpern method for cocoercive operators and their
two-step extrapolated Halpern method for general Lipschitz-monotone operators
both attain \(O(\epsilon^{-3})\) oracle complexity for an expected residual;
under sharpness, restarting gives
\(O(\epsilon^{-2}\log(1/\epsilon))\).

On gradient fields, the Baillon--Haddad theorem identifies
\(1/L\)-cocoercivity with convex
\(L\)-smoothness~\cite{baillon1977,bauschke2010baillon}. Cocoercivity alone does
not imply the existence of a potential~\cite{rockafellar1966characterization}.
Even on the common subclass, their
operator residual and our objective gap are not equivalent at the same order:
\[
        \E\|\nabla f(X_N)\|
        \le
        \sqrt{2L\,\E[f(X_N)-f(x^*)]},
\]
and the converse needs an error bound. We therefore make no improvement claim
for their oracle complexity. Our result instead uses one oracle estimate per
iteration, permits an \(\F_n\)-measurable bias and growing raw variance under
trajectory-wise scaled control, and gives strong selection in Hilbert space;
their stronger oracle yields sharper complexity for a broader operator class.

\subsection{Other stochastic fixed-point frameworks}
\label{subsec:other-stochastic-fixed-point}

Bravo--Cominetti~\cite{bravocominetti2024} obtain a.s. convergence to a random
fixed point and expected-residual bounds for perturbed Krasnosel'ski\u{\i}--Mann
iterations, allowing suitably controlled growing variance.
Pischke--Powell~\cite{pischke2024asymptotic} give quantitative asymptotic
regularity and oracle-complexity results for generalized stochastic
Halpern--Mann schemes. These fixed-operator frameworks are broader; our smooth
gradient structure instead yields strong convergence to the metric projection
of the anchor and objective last-iterate bounds.

\subsection{Summary of assumptions and conclusions}
\label{subsec:comparison-table}

\begin{table}[htbp]
\centering
\scriptsize
\renewcommand{\arraystretch}{1.12}
\begin{tabular}{@{}
>{\raggedright\arraybackslash}p{0.17\textwidth}
>{\raggedright\arraybackslash}p{0.20\textwidth}
>{\raggedright\arraybackslash}p{0.26\textwidth}
>{\raggedright\arraybackslash}p{0.23\textwidth}@{}}
\toprule
Reference / framework
& Problem setting
& Noise / oracle assumptions
& Main conclusion
\\
\midrule
Alacaoglu--Malitsky--Wright~\cite{alacaoglu2025weakervar}
& Constrained, possibly nonsmooth, convex minimization in \(\mathbb R^d\)
& Unbiased subgradient oracle with a Blum--Gladyshev second-moment bound
& \(O(k^{-1/2})\) weighted-average and
  \(\widetilde O(k^{-1/2})\) last-iterate rates  in expectation 
\\
\midrule
Bravo--Contreras~\cite{bravo2024stochastic}
& Stochastic Halpern iteration for fixed nonexpansive or contractive maps
& Unbiased operator oracle and increasing minibatches
& \(\widetilde O(N^{-1})\) expected residual,
  \(\widetilde O(\epsilon^{-5})\) oracle complexity, and a.s. convergence under
  summable mean errors
\\
\midrule
Cai--Song--Guzm\'an--Diakonikolas~\cite{cai2022stochastic}
& Finite-dimensional cocoercive or Lipschitz-monotone equations
& Unbiased bounded-variance multipoint oracle, mean-square Lipschitz condition,
  and PAGE
& \(O(\epsilon^{-3})\) expected-residual oracle complexity;
  \(O(\epsilon^{-2}\log(1/\epsilon))\) under sharpness
\\
\midrule
Bravo--Cominetti~\cite{bravocominetti2024}
& Perturbed Krasnosel'ski\u{\i}--Mann iterations for fixed nonexpansive maps
& Martingale-difference perturbations with weighted variance control
& A.s. convergence to a random fixed point and expected-residual bounds
\\
\midrule
Pischke--Powell~\cite{pischke2024asymptotic}
& Generalized stochastic Halpern--Mann schemes
& Abstract stochastic errors with quantitative decay assumptions
& Rates of asymptotic regularity and oracle-complexity results
\\
\midrule
This paper: pathwise regime
& Smooth convex minimization with selected solution \(P_S(u)\)
& One oracle estimate per step; martingale noise and an
  \(\F_n\)-measurable bias satisfying \eqref{eq:comparison-summability}
& A.s. strong convergence to \(P_S(u)\)
\\
\midrule
This paper: rate regime
& The same smooth convex problem
& 
Global and suffix controls; additionally, in the unbiased critical regime,
  a uniform conditional \(p\)-moment bound with \(p>4\)
& 
\(\widetilde O(k^{-1/2})\) expected last-iterate rate; under the additional
  \(p\)-moment condition, also a.s.\ strong convergence
\\
\bottomrule
\end{tabular}
\caption{Schematic comparison with related stochastic Halpern and fixed-point
frameworks.}
\label{tab:comparison-literature}
\end{table}
\FloatBarrier

\section{Conclusion}
\label{sec:conclusion}

Halpern anchoring selects \(P_S(u)\), while the vanishing gradient step scales
the \(\F_n\)-measurable bias and martingale fluctuation. Pathwise comparison
summability yields a.s. strong convergence; expectation-level global and
suffix control yields last-iterate bounds. 
In the unbiased critical regime
\(\varepsilon_n=n^{-1/2}\),
\(\alpha_n=(\log n)/n\) eventually, a uniform conditional
\((4+\delta)\)-moment bound combines the two conclusions.

The counterexamples below
separate these mechanisms but do not claim necessity of our precise sufficient
conditions.

Open questions include removing logarithmic losses and deriving suffix control
from structural oracle hypotheses. The threshold in
\cref{thm:moment-as-and-rate} may also be improvable: its proof uses \(p>4\),
whereas any uniform moment of order greater than two already excludes the
fixed-size jumps in \cref{ex:role-anti-spike}.

\appendix

\section{Auxiliary summability lemmas}\label{app:lemmas}

This appendix collects the auxiliary tools used in the paper: a tail-wise
weighted telescoping lemma, a suffix-averaging lemma, and the harmonic estimates
needed to verify the suffix-regularity conditions.

\begin{lemma}
\label[lemma]{lem:tail-weighted-telescoping}
Let $N\le M$ be integers and let
\((a_n)_{n=N}^{M+1}\) and \((r_n)_{n=N}^{M}\) be nonnegative finite real
sequences, and let \((b_n)_{n=N}^{M}\) be a finite real sequence.
Let $w_n>0$ and $\varrho_n>0$ for
$N\le n\le M$, and assume that
\(a_{n+1}+w_n b_n\le\varrho_n a_n+r_n, \qquad N\le n\le M.\)
Define the factors $(\Lambda_n)_{n=N-1}^{M}$ by
$\Lambda_{N-1}:=1,$ $\Lambda_n\varrho_n=\Lambda_{n-1},$ $N\le n\le M,$
equivalently $\Lambda_n=\prod_{k=N}^{n}\varrho_k^{-1}$. Then
\(\sum_{n=N}^{M}\Lambda_n w_n b_n\le a_N+\sum_{n=N}^{M}\Lambda_n r_n,\)
and in particular
\[
  \min_{N\le n\le M} b_n
  \le\frac{a_N+\sum_{n=N}^{M}\Lambda_n r_n}{\sum_{n=N}^{M}\Lambda_n w_n}.
\]
Moreover, let $(a_n)$ and $(r_n)$ be nonnegative real sequences, let
$(b_n)$ be a real sequence, and let $(w_n)$ be a positive
sequence satisfying
\(a_{n+1}+w_n b_n\le a_n+r_n.\)
Let $M(N)\ge N$ satisfy $\sum_{n=N}^{M(N)}w_n\to+\infty$. If
$\sup_{n\ge0}a_n<+\infty$ and
\(\frac{\sum_{n=N}^{M(N)}r_n}{\sum_{n=N}^{M(N)}w_n}\longrightarrow0,\)
then \(\limsup_{N\to\infty} \min_{N\le n\le M(N)} b_n\le0.\)
If, in addition, \(b_n\ge0\) for all \(n\), then
\(\min_{N\le n\le M(N)}b_n\to0\).
\end{lemma}

\begin{proof}
Multiplication by \(\Lambda_n\) and summation give
\[
 \sum_{n=N}^M\Lambda_nw_nb_n
 \le a_N-\Lambda_Ma_{M+1}+\sum_{n=N}^M\Lambda_nr_n
 \le a_N+\sum_{n=N}^M\Lambda_nr_n.
\]
The minimum bound follows because all \(\Lambda_nw_n\) are positive. For the
asymptotic claim set \(\varrho_n=1\): the bound becomes
\(\min_{N\le n\le M(N)}b_n\le a_N/W_N+R_N/W_N\), where
\(W_N=\sum w_n\to\infty\) and \(R_N/W_N\to0\). Nonnegativity gives the final
assertion.
\end{proof}

We turn to the suffix-averaging lemma. It is inspired by
last-iterate and suffix-averaging arguments in
Orabona~\cite{orabona2020last},
Shamir--Zhang~\cite{shamir2013stochastic}, and
Alacaoglu--Malitsky--Wright~\cite{alacaoglu2025weakervar}.

\begin{lemma}
\label[lemma]{lem:suffix-averaging}
Let $(d_n)_{n\ge0}$ be a real sequence. For every \(K\ge1\),
\begin{equation}
\label{eq:suffix-averaging-identity}
  d_K
  =
  \frac1K\sum_{k=1}^{K}d_k
  +
  \sum_{\ell=1}^{K-1}
  \frac1{\ell(\ell+1)}
  \sum_{k=K-\ell}^{K}
  \bigl(d_k-d_{K-\ell}\bigr).
\end{equation}
If, in addition, $(q_n)_{n\ge0}$ is nonnegative and
\[
        [d_{n+1}-d_n]_+\le q_n,
        \qquad n\ge0,
        \qquad [r]_+:=\max\{r,0\},
\]
then, for every $0\le m\le K$,
\[
  d_K
  \le\frac1{K-m+1}\sum_{n=m}^{K}d_n
  +\frac1{K-m+1}\sum_{n=m}^{K-1}(n-m+1)q_n.
\]
Consequently,
\[
  d_K
  \le\inf_{0\le m\le K}
    \left\{\frac1{K-m+1}\sum_{n=m}^{K}d_n
      +\frac1{K-m+1}\sum_{n=m}^{K-1}(n-m+1)q_n\right\}.
\]
\end{lemma}

\begin{proof}
The coefficients of \(d_K\) and \(d_j\), \(j<K\), on the right of
\eqref{eq:suffix-averaging-identity} are respectively
\(K^{-1}+\sum_{\ell<K}[\ell(\ell+1)]^{-1}=1\) and zero, by telescoping.
Moreover, \(d_K\le d_j+\sum_{n=j}^{K-1}q_n\) for \(m\le j\le K\). Averaging
over \(j\) and interchanging the finite sums gives
\[
 d_K\le\frac{\sum_{j=m}^{K}d_j+
 \sum_{n=m}^{K-1}(n-m+1)q_n}{K-m+1}.
\]
Take the infimum over \(m\).
\end{proof}

We also record the harmonic estimates used for suffix regularity.

\begin{lemma}
\label[lemma]{lem:polylog-sums}
Let $s\ge0$. Then, as $K\to+\infty$,
\(\sum_{n=1}^{K}\frac{\log^s(n+1)}{n}=O\bigl(\log^{s+1}(K+1)\bigr),\)
and, for $0<p<1$,
\(\sum_{n=1}^{K}\frac{\log^s(n+1)}{n^p}=O\bigl(K^{1-p}\log^s(K+1)\bigr).\)
Moreover, for $1\le\ell\le K$,
\[
  \sum_{n=K-\ell+1}^{K}\frac{\log^s(n+1)}{n}
  \le C_s
  \begin{cases}
    \dfrac{\ell\,\log^s(K+1)}{K}, & 1\le\ell\le K/2,\\[1.2em]
    \log^{s+1}(K+1), & K/2<\ell\le K,
  \end{cases}
\]
for a constant $C_s>0$ independent of $K$ and $\ell$.
Consequently,
\begin{equation}
\label{eq:polylog-suffix-double-sum}
  \sum_{\ell=1}^{K-1}\frac1{\ell(\ell+1)}
  \sum_{n=K-\ell}^{K}\frac{\log^s(n+1)}{n}
  =O\!\left(\frac{\log^{s+1}(K+1)}{K}\right).
\end{equation}
\end{lemma}

\begin{proof}
The global estimates follow from the integral test. For the local one, use
\(n\ge K/2\) if \(\ell\le K/2\), and the global bound otherwise. Splitting the
outer sum in \eqref{eq:polylog-suffix-double-sum} at \(K/2\) then gives a
harmonic sum divided by \(K\) on the first part and an \(O(K^{-1})\) tail on
the second.
\end{proof}

\section{Counterexamples}
\label{app:counterexamples}

All examples use
\[
        \HH=\mathbb R, \qquad f(x)=\frac12x^2, \qquad
        S=\{0\}, \qquad u=0, \qquad x^*=0.
\]
Then \(\nabla f(x)=x\), and the stochastic Halpern--gradient recursion becomes
\begin{equation}
\label{eq:scalar-quadratic-recursion}
        X_{n+1}
        =
        (1-\alpha_n)
        \bigl(
            (1-\varepsilon_n)X_n-\varepsilon_n U_n
        \bigr).
\end{equation}
Here \(X_0\) is deterministic and the coefficients satisfy the standing
conditions unless specified otherwise.

\begin{example}
\label[example]{ex:rare-large-spikes}
Let \((m_j)\) be a strictly increasing sequence of positive integers and let
\((s_j)\) be a positive sequence with \(s_j\to0\). Let
\((\zeta_j)\) be independent Rademacher random variables, that is,
\(\mathbb P(\zeta_j=1)=\mathbb P(\zeta_j=-1)=\frac12.\)
Equip the construction with the completed filtration
\[
        \F_n
        :=
        \sigma\bigl(X_0,\zeta_j:m_j\le n\bigr),
        \qquad n\ge0.
\]
Define the noise by
\begin{equation}
\label{eq:sparse-rademacher-noise}
        U_{m_j-1}
        =
        \frac{\sqrt{s_j}}{\varepsilon_{m_j-1}}\zeta_j,
        \qquad
        U_n=0
        \quad\text{if } n\neq m_j-1.
\end{equation}
Then
\(U_n\) is \(\F_{n+1}\)-measurable, belongs to \(L^2\) for every
fixed \(n\), and \(\E[U_n\mid\F_n]=0\). The recursion therefore
defines an adapted \(L^2\) process by induction.
At a spike time,
\[
\begin{aligned}
        v_{m_j-1}
        &:=
        \varepsilon_{m_j-1}^2
        \E[
            |U_{m_j-1}|^2
            \mid\F_{m_j-1}
        ]                                                        \\
        &=
        \varepsilon_{m_j-1}^2
        \frac{s_j}{\varepsilon_{m_j-1}^2}
        \E|\zeta_j|^2
        =
        s_j.
\end{aligned}
\]
Thus the weighted variance contribution of the \(j\)-th spike is exactly
\(s_j\).

At the corresponding iterate \(m_j\), the recursion gives
\[
        X_{m_j}
        =
        (1-\alpha_{m_j-1})
        \bigl(
            (1-\varepsilon_{m_j-1})X_{m_j-1}
            -
            \sqrt{s_j}\,\zeta_j
        \bigr).
\]
Since \(\zeta_j\) is independent of \(\F_{m_j-1}\) and centered,
\[
\begin{aligned}
        \E[
            X_{m_j}^2
            \mid
            \F_{m_j-1}
        ]
        &=
        (1-\alpha_{m_j-1})^2
        \left[
            (1-\varepsilon_{m_j-1})^2X_{m_j-1}^2+s_j
        \right]                                                   \\
        &\ge
        (1-\alpha_{m_j-1})^2s_j.
\end{aligned}
\]
Choosing \(j\) large enough so that \(\alpha_{m_j-1}\le1/2\), we obtain
\(\E X_{m_j}^2 \ge \frac14s_j.\)
Consequently,
\begin{equation}
\label{eq:spike-gap-lower-bound}
        \E[f(X_{m_j})-f(x^*)]
        =
        \frac12\E X_{m_j}^2
        \ge
        \frac18s_j
\end{equation}
for all sufficiently large \(j\).

Thus a single spike can be calibrated exactly: choosing the effective
displacement to have squared size \(s_j\) makes the weighted variance
contribution equal to \(s_j\), and it also forces
\(\E[f(X_{m_j})-f(x^*)]\ge c s_j\)
at the next iterate. The later examples use this calibration to place errors
at prescribed terminal indices while keeping the total weighted variance under
control.
\end{example}

\begin{example}
\label[example]{ex:WS-no-last-iterate-rate}
This example shows exactly what weighted summability does not control. Given
any prescribed positive rate scale \((\varphi_n)\) with \(\varphi_n\downarrow0\), we
construct an unbiased sparse-noise sequence for which
\(\sum_{n=0}^{\infty}v_n<+\infty,\)
but the expected terminal gaps fail to satisfy
\(\E[f(X_n)-f(x^*)]=O(\varphi_n).\)
Thus even summable weighted variance does not preclude rare terminal iterates
whose error is much larger than a prescribed last-iterate scale.

Let \((\varphi_n)\) be any positive sequence with \(\varphi_n\downarrow0\). We show
that one can choose spike locations \((m_j)\) and spike sizes \((s_j)\) such
that
\(\sum_{j=1}^{\infty}s_j<+\infty,\)
but
\(\frac{s_j}{\varphi_{m_j}}\to+\infty.\)
Indeed, choose \(m_j\) so large that
\(\varphi_{m_j}\le \frac1{j^3},\)
and set
\(s_j:=j\varphi_{m_j}.\)
Then
\(s_j\le \frac1{j^2},\)
so \(\sum_js_j<+\infty\), while
\(\frac{s_j}{\varphi_{m_j}}=j\to+\infty.\)

Now define the noise by \eqref{eq:sparse-rademacher-noise}, with the completed
filtration specified in \cref{ex:rare-large-spikes}. Since
\(v_{m_j-1}=s_j, \qquad v_n=0 \quad\text{otherwise},\)
we have
\begin{equation}
\label{eq:WS-counterexample}
        \sum_{n=0}^{\infty}v_n
        =
        \sum_{j=1}^{\infty}s_j
        <+\infty.
\end{equation}
Thus the weighted summability condition holds.

However, by \eqref{eq:spike-gap-lower-bound},
\(\E[f(X_{m_j})-f(x^*)] \ge c s_j\)
for some \(c>0\), and therefore
\[
        \frac{
            \E[f(X_{m_j})-f(x^*)]
        }{
            \varphi_{m_j}
        }
        \ge
        c\frac{s_j}{\varphi_{m_j}}
        \to+\infty.
\]
Hence the estimate
\(\E[f(X_n)-f(x^*)]=O(\varphi_n)\)
fails.

If one only wants to rule out a little-\(o\) estimate, choose instead
\(s_j:=\varphi_{m_j},\)
with the indices \(m_j\) so sparse that
\(\sum_j\varphi_{m_j}<+\infty.\)
Then weighted summability still holds, while
\(\E[f(X_{m_j})-f(x^*)] \ge c\varphi_{m_j}.\)
Thus
\(\E[f(X_n)-f(x^*)]=o(\varphi_n)\)
fails.

Thus the same construction can be tuned either to violate a prescribed
\(O(\varphi_n)\) bound or, with smaller spikes, to violate a prescribed
\(o(\varphi_n)\) bound. The point is not terminology: the summability condition
controls the total variance budget \(\sum_n v_n\), but it does not control how
that budget is distributed along the terminal indices at which a last-iterate
rate is tested.
\end{example}

\begin{example}
\label[example]{ex:suffix-no-as-convergence}
This example constructs an unbiased bounded-variance noise sequence for which
the stochastic suffix-regularity estimates hold, and hence the last-iterate
bound in expectation is available, but the iterates do not converge almost
surely. The construction uses infinitely many jumps of fixed effective size
\(\varepsilon_n|U_n|=a\), occurring with probabilities \(1/(n+1)\).

Consider the logarithmic-anchor regime
\(\varepsilon_n=(n+1)^{-1/2}\) and
\(\alpha_n=\frac{\log(n+2)}{n+1}\).
Let \((B_n)\) be independent Bernoulli random variables with
\(\mathbb P(B_n=1)=\frac1{n+1},\) \( \mathbb P(B_n=0)=1-\frac1{n+1},\)
and let \((\zeta_n)\) be independent Rademacher random variables, independent
of \((B_n)\). Assume that all these variables are independent of \(X_0\), and
use the completed filtration
\[
        \F_n
        :=
        \sigma\bigl(X_0,(B_k,\zeta_k):0\le k<n\bigr).
\]
Thus \((B_n,\zeta_n)\) is independent of \(\F_n\) and is revealed
between times \(n\) and \(n+1\). Fix \(a>0\), and define
\begin{equation}
\label{eq:bernoulli-spike-noise}
        U_n
        =
        \frac{a}{\varepsilon_n}\zeta_n B_n.
\end{equation}
Then
\(U_n\) is \(\F_{n+1}\)-measurable and
\(\E[U_n\mid\F_n]=0\).
Moreover,
\(
        \E[
            |U_n|^2
            \mid\F_n
        ]
        =
        \frac{a^2}{\varepsilon_n^2}
        \En[B_n^2]                                            =
        \frac{a^2}{\varepsilon_n^2}
        \frac1{n+1}.
\)
Since \(\varepsilon_n^2=(n+1)^{-1}\), we get
\(\E[ |U_n|^2 \mid\F_n ] = a^2.\)
Thus the stochastic-gradient noise has uniformly bounded conditional
variance. Moreover,
\(v_n = \varepsilon_n^2 \E[ |U_n|^2 \mid\F_n ] = \frac{a^2}{n+1}.\)
In particular, \(\E v_n\le C\alpha_n\); since the bias is zero,
\cref{lem:L2-boundedness-stoch} gives
\(\sup_n\E|X_n|^2<+\infty\). Hence all the integrability hypotheses of
the last-iterate theorem are satisfied.
Therefore, in the notation of the stochastic last-iterate section,
\[
        \theta_n
        =
        \alpha_n+\varepsilon_n^2+\E v_n
        =
        O\left(\frac{\log(n+2)}{n+1}\right).
\]
Consequently, \cref{lem:polylog-sums} with \(s=1\) gives
\[
        \frac{
            1+\sum_{k=1}^{K}\theta_k
        }{K}
        =
        O\left(\frac{(\log K)^2}{K}\right), \text{ and} \quad
        \sum_{\ell=1}^{K-1}
        \frac1{\ell(\ell+1)}
        \sum_{k=K-\ell}^{K}\theta_k
        =
        O\left(\frac{(\log K)^2}{K}\right).
\]
Since \(\varepsilon_K=(K+1)^{-1/2}\asymp K^{-1/2}\), the
suffix-regularity condition holds.
In particular, \cref{thm:stoch-last-iterate-general} gives the corresponding
last-iterate estimate in expectation for this construction.

However,
\(\sum_{n=0}^{\infty}\mathbb P(B_n=1) = \sum_{n=0}^{\infty}\frac1{n+1} = +\infty.\)
By the second Borel--Cantelli lemma, since the events \(\{B_n=1\}\) are
independent,
\(B_n=1 \) for infinitely many \(n\) a.s.
At every such spike time, the recursion reads
\(X_{n+1} = (1-\alpha_n) \bigl( (1-\varepsilon_n)X_n-a\zeta_n \bigr).\)
Suppose, on the contrary, that \((X_n)\) converges to a finite random limit.
Then \((X_n)\) is bounded and \(X_{n+1}-X_n\to0\). At a spike time,
\[
\begin{aligned}
        X_{n+1}-X_n
        &=-\bigl(\alpha_n+\varepsilon_n-\alpha_n\varepsilon_n\bigr)X_n
          -(1-\alpha_n)a\zeta_n.
\end{aligned}
\]
Along the infinitely many spike times, the first term tends to zero because
\(\alpha_n,\varepsilon_n\to0\) and \((X_n)\) is bounded, whereas the modulus
of the second term tends to \(a\). Hence
\(|X_{n+1}-X_n|\to a\) along those times, a contradiction. Therefore the
iterates do not converge almost surely.

Thus this noise satisfies the expectation-level suffix estimates used for
last-iterate bounds, while the actual trajectory does not converge. The example
separates the expectation rate mechanism from the pathwise mechanism needed
for almost sure convergence.
\end{example}

\begin{example}
\label[example]{ex:role-anti-spike}
This example identifies which part of the preceding Bernoulli construction is
excluded by a uniform \(4+\delta\) moment bound. The noise in
\eqref{eq:bernoulli-spike-noise} has bounded conditional second moment, but
its conditional \(4+\delta\) moments diverge. Indeed, for every \(\delta>0\),
\[
        \E[
            |U_n|^{4+\delta}
            \mid\F_n
        ]
        =
        \left(\frac{a}{\varepsilon_n}\right)^{4+\delta}
        \En[B_n]                                              =
        a^{4+\delta}
        \varepsilon_n^{-(4+\delta)}
        \frac1{n+1}.
\]
Since
\(\varepsilon_n=(n+1)^{-1/2},\)
this becomes
\[
        \E[
            |U_n|^{4+\delta}
            \mid\F_n
        ]
        =
        a^{4+\delta}
        (n+1)^{2+\delta/2}
        \frac1{n+1}
        =
        a^{4+\delta}(n+1)^{1+\delta/2}
        \to+\infty.
\]
Therefore the spike construction satisfies bounded variance but violates
every uniform \(4+\delta\) moment bound.

Conversely, a uniform \(4+\delta\) moment condition prevents infinitely many
large effective jumps of the form
\(\varepsilon_n|U_n|\ge a.\)
Indeed, if
\(\sup_n \E[ |U_n|^{4+\delta} \mid\F_n ] \le C \) a.s.,
then Markov's inequality gives
\[
        \mathbb P(
            \varepsilon_n|U_n|\ge a
            \mid \F_n
        )
        =
        \mathbb P(
            |U_n|\ge a/\varepsilon_n
            \mid \F_n
        )
        \le
        \frac{
            \varepsilon_n^{4+\delta}
        }{
            a^{4+\delta}
        }
        \E[
            |U_n|^{4+\delta}
            \mid\F_n
        ]
        \le
        \frac{C}{a^{4+\delta}}\varepsilon_n^{4+\delta}.
\]
In the logarithmic-anchor regime
\(\varepsilon_n=(n+1)^{-1/2},\)
one has
\[
        \sum_{n=0}^{\infty}\varepsilon_n^{4+\delta}
        =
        \sum_{n=0}^{\infty}(n+1)^{-2-\delta/2}
        <+\infty.
\]
Taking expectations and summing the preceding bound shows that the
probabilities of the large-effective-jump events are summable. Hence the first
Borel--Cantelli lemma excludes infinitely many such events. The same
calculation works for any moment exponent \(r>2\), since
\[
        \sum_{n=0}^{\infty}\varepsilon_n^r
        =\sum_{n=0}^{\infty}(n+1)^{-r/2}<+\infty.
\]

Thus a uniform \(4+\delta\) moment bound excludes the particular fixed-size
effective jumps used in \cref{ex:suffix-no-as-convergence}; in fact, any
uniform moment of order strictly larger than two excludes those individual
jumps. The threshold \(p>4\) in \cref{thm:moment-as-and-rate} is additionally
used to control accumulated fluctuations over logarithmic blocks. In this
sense the higher-moment hypothesis provides an anti-spike control stronger
than bounded variance, but the example does not prove that the precise
threshold \(p>4\) is necessary.
\end{example}

\section*{Declarations}

The authors used a large language model to assist in grammar checking the manuscript
and in condensing its presentation. The authors verified all mathematical
statements, proofs, citations, and final wording. The authors assume
responsibility for all content.

\bibliographystyle{siamplain}
\bibliography{references}

\end{document}